\input ./style/arxiv-general.cfg
\documentclass[aos,MSNbibl,seceqn,dvips]{arximspdf}
\makeatletter
   \@ifpackageloaded{graphicx}{}{\usepackage{graphicx}}
\makeatother

\doi{10.1214/15-AOS1341}
\volume{43}
\issue{5}
\pubyear{2015}
\firstpage{2259}
\lastpage{2295}
\docsubty{FLA}

\makeatletter

\def\vfrac#1#2{(#1)/#2}

\def\sklfrac#1#2{(#1/#2)}

\newcommand{\rrvert}{\vert}
\newcommand{\rrVert}{\Vert}
\newcommand{\llvert}{\vert}
\newcommand{\llVert}{\Vert}
\renewcommand{\mid}{|}
\newcommand{\N}{\mathbb{N}}
\newcommand{\Z}{\mathbb{Z}}
\newcommand{\utheta}{\underline{\theta}}
\newcommand{\uY}{\underline{Y}}
\newcommand{\R}{\mathbb{R}}
\newtheorem{teo}{Theorem}[section]
\newproclaim{defi}{Definition}
\newtheorem{prop}{Proposition}
\newtheorem{lemma}{Lemma}
\newproclaim{rem}{Remark}
\makeatother

\begin{document}
\begin{frontmatter}

\title{On adaptive posterior concentration rates}
\runtitle{On adaptive posterior concentration rates}

\begin{aug}
\author[A]{\fnms{Marc}~\snm{Hoffmann}\corref{}\thanksref{M1,T1}\ead[label=e1]{hoffmann@ceremade.dauphine.fr}},
\author[A]{\fnms{Judith}~\snm{Rousseau}\thanksref{M1,T1}\ead[label=e2]{rousseau@ceremade.dauphine.fr}}\\
\and
\author[B]{\fnms{Johannes}~\snm{Schmidt-Hieber}\thanksref{M2,T3}\ead[label=e3]{schmidthieberaj@math.leidenuniv.nl}}
\runauthor{M. Hoffmann, J. Rousseau and J. Schmidt-Hieber}
\affiliation{Universit\'e Paris-Dauphine\thanksmark{M1} and
University of Leiden \thanksmark{M2}}
\address[A]{M. Hoffmann\\
J. Rousseau\\
CEREMADE\\
Universit\'e Paris Dauphine\\
Place du Mar\'echal De Lattre de Tassigny\\
75016 Paris\\
France\\
\printead{e1}\\
\phantom{E-mail:\ }\printead*{e2}}
\address[B]{Mathematical Institute\\
University of Leiden\\
Niels Bohrweg 1\\
2333 CA Leiden\\
The Netherlands\\
\printead{e3}}
\end{aug}
\thankstext{T1}{Supported in part by the Agence Nationale de la Recherche
(project BANDHIT and Blanc SIMI 1~2011 project CALIBRATION).}
\thankstext{T3}{Supported by DFG postdoctoral fellowship SCHM 2807/1-1.}

%
\received{\smonth{5} \syear{2013}}
%
\revised{\smonth{4} \syear{2015}}

%
\begin{abstract}
We investigate the problem of deriving posterior concentration rates
under different loss functions in nonparametric Bayes. We first provide
a lower bound on posterior coverages of shrinking neighbourhoods that
relates the metric or loss under which the shrinking neighbourhood is
considered, and an intrinsic pre-metric linked to frequentist
separation rates.
In the Gaussian white noise model, we construct feasible priors based
on a spike and slab procedure reminiscent of wavelet thresholding that
achieve adaptive rates of contraction under $L^2$ or $L^\infty$
metrics when the underlying parameter belongs to a collection of H\"
older balls and that moreover achieve our lower bound. We analyse the
consequences in terms of asymptotic behaviour of posterior credible
balls as well as frequentist minimax adaptive estimation. Our results
are appended with an upper bound for the contraction rate under an
arbitrary loss in a generic regular experiment. The upper bound is
attained for certain sieve priors and enables to extend our results to
density estimation.
\end{abstract}

%
\end{frontmatter}

\setcounter{footnote}{2}

\section{Introduction}
\subsection{Setting}
There has been a growing interest for posterior concentration rates in
nonparametric Bayes over the last decade, initiated by the seminal
papers of Schwartz \cite{schwartz}, Barron \cite{barron1988} and
Ghosal, Ghosh and van~der Vaart
\cite{ghosalghoshvdv00}.
Consider a statistical model or experiment $\mathcal E^n = \{P_\theta
^n:\theta\in\Theta\}$ generated by data $Y^n$, with parameter space
$\Theta$ equipped with a prior distribution~$\pi$.
The posterior distribution $P^\pi(\cdot\mid  Y^n)$ concentrates at rate
$\epsilon_n>0$ under $P_{\theta_0}^n$ for the loss $\ell: \Theta
\times\Theta\rightarrow[0,\infty)$ if
%
\begin{equation}
\label{defpostconc} E_{\theta_0}^n \bigl[ P^\pi \bigl(
\theta: \ell( \theta, \theta_0) >\epsilon_n \mid
Y^n \bigr) \bigr] = o(1).
\end{equation}
%
Posterior concentration allows to uncover frequentist properties of
Bayesian methods. It implies that the posterior probability of an
$\epsilon_n$-neighbour\-hood around the true parameter $\theta_0$
converges to one. Thus, most of the posterior mass will be close to the
truth in the frequentist sense.

Whenever (\ref{defpostconc}) holds uniformly in $\theta_0 \in
\Theta$ and if $\epsilon_n$ can be taken as constant multiple of the
minimax rate of estimation over $\Theta$ for the loss $\ell$, we say
that the concentration rate is asymptotically minimax.
We further say that the posterior distribution $P^\pi(\cdot\mid Y^n)$
concentrates adaptively over the collection $\{\Theta_\beta, \beta
\in\mathcal I\}$ of subsets of $\Theta$ if
%
\begin{equation}
\label{defpostadap} \sup_{\theta_0 \in\Theta_\beta} E_{\theta_0}^n
\bigl[ P^\pi \bigl(\theta: \ell( \theta, \theta_0) >
\epsilon_n(\beta) \mid Y^n \bigr) \bigr] = o(1)\qquad
\mbox{for every }\beta\in{\mathcal I},
\end{equation}
where $\epsilon_n(\beta)$ is a constant multiple of the minimax rate
of adaptive estimation over $\Theta_\beta$. Recently, some families
of prior distributions under various types of statistical models have
been studied in this light and have been proved to lead to adaptive
posterior concentration rates; see Section~\ref{furtherbiblio} for a
more extensive discussion of these results.
Similarly as for (\ref{defpostconc}), existence of a result of type
(\ref{defpostadap}) implies that the Bayes estimator is minimax
adaptive under fairly general conditions; see Section~\ref
{freq-consequences}. Consequently, existence or nonexistence of
adaptive estimators in some nonparametric situations (see, e.g.,
\cite{lepski90,brownlow96}) yield limitations about
the best possible achievable concentration rates $\epsilon_n(\beta)$
in (\ref{defpostadap}).

In this paper, we are interested in understanding further the interplay
between nonparametric minimax rates of convergence and the existence of
adaptive concentration rates for appropriate priors in nonparametric estimation.
We cover in particular the two paradigmatic examples of density
estimation, when the data $Y^n$ is drawn from a $n$-sample of an
unknown distribution, and the case of a signal observed in Gaussian
white noise.
More specifically, we attempt to answer the following related questions:
\begin{longlist}[(II)]
\item[(I)] Can we formalise the connexion between posterior
concentration rates and the minimax theory: given an experiment
$\mathcal E^n$ and a loss $\ell$, can we define some notion of lower
bound associated to the posterior concentration rate? Can we derive a
generic construction for a prior with posterior achieving the minimax
rate of convergence in the sense of (\ref{defpostconc})? Can we
further make this construction adaptive, in the sense of (\ref{defpostadap})?
\item[(II)] In the specific framework of the $L^2$ and $L^\infty$
metric for the loss $\ell$, can we construct a feasible prior in
standard models such as Gaussian white noise or density estimation for
which the posterior distribution contracts adaptively over H\"older balls?
\end{longlist}


\subsection{Main results} \label{mainresults}
A first answer to these problems is given in Section~\ref{seclower}
in the form of a lower bound on the speed at which the posterior mass
outside an $\epsilon_n$-ball vanishes in the sense of (\ref
{defpostconc}). Assume that $\Theta$ is equipped with a
pre-metric\footnote{That is, we only require that $d$ is
nonnegative and $d(\theta, \theta')=0$ iff $\theta= \theta'$.} $d$
that controls the separation rate between two elements in $\mathcal
E^n$. We prove in Theorem~\ref{thlowerbound} that if $\mathcal E^n$
is dominated and admits a certain regularity condition
then, for every prior $\pi$ such that the posterior $P^\pi(\cdot
\mid Y^n)$ concentrates with rate $\epsilon_n$ over $\Theta$, there
exists a constant $c>0$ such that
%
\begin{equation}
\label{lowerboundinformal} \sup_{\theta_0 \in\Theta}   E_{\theta_0}^n
\bigl[ P^\pi \bigl(\theta: \ell(\theta, \theta_0 ) \geq
\epsilon_n \mid Y^n \bigr) \bigr] \geq e^{-c n\Omega(\epsilon_n, \Theta, \ell)^2},
\end{equation}
where
\[
\Omega(\epsilon_n,\Theta,\ell)=\inf\bigl\{d\bigl(\theta,
\theta'\bigr):  \ell\bigl(\theta, \theta'\bigr) \geq2
\epsilon_n,\theta, \theta' \in\Theta \bigr\}.
\]
The pre-metric $d$ geometrises the statistical model and does not
depend on the loss function nor the rate. At this point, one might
think of $d$ as the Hellinger distance. If $\epsilon_n \rightarrow0$,
only the local behaviour of $d$ plays a role in the definition of
$\Omega(\epsilon_n, \Theta, \ell)$, which gives slightly more
flexibility and allows to take, for instance, $d$ as the $L^2$-metric in
the Gaussian white noise model. The precise conditions that determine
$d$ are stated in Theorems~\ref{thlowerbound}~and~\ref{thuppbound}. Explicit computations are developed in Section~\ref{seclower}.

The exponent $\Omega(\epsilon_n, \Theta, \ell)$ appearing in (\ref
{lowerboundinformal}) is a dual formulation of the modulus of
continuity introduced in \cite{donoholiu} and further considered by
Cai and Low \cite{cailow06},
Cai, Low and Zhao \cite{CaiLowZhao04}; see Section~\ref{furtherbiblio}.
Theorem~\ref{thlowerbound} also admits a stronger local version:
$\epsilon_n$ can be a function of $\theta$ also in a manner similar
to the between classes modulus of continuity of
Cai and Low \cite{cailow04}.
Another important consequence is that there are limitations of the
commonly employed proof strategy for derivation of posterior
concentration rates; see Section~\ref{subsecconseq}.

In Section~\ref{secWN}, we address question (II) and explicitly
construct a prior---in the family of spike and slab priors---that
achieves the lower bound of Theorem~\ref{thlowerbound} in the white
noise model simultaneously over a collection of H\"older balls
$\mathcal H(\beta, L)$ for $\beta\in{\mathcal I}$, where ${\mathcal
I}$ is a compact subset of $(0,\infty)$. Recasting $Y^n$ into a
regular wavelet basis (see, e.g., \cite{cohen2000} and \cite
{cdv1998}), we obtain the sequence model
\[
Y_{j,k}=\theta_{j,k}+n^{-1/2}\epsilon_{j,k},\qquad
  k\in I_j,j=0,1,\ldots,
\]
where $k \in I_j$ is a location parameter at scale $2^{-j}$ with $I_j$
having approximately $2^j$ terms, and the $\epsilon_{j,k}$ are i.i.d.
standard normal.
The spike and slab prior is constructed as follows: for $j$ less than a
maximal resolution level\footnote{In the sequel, we adopt the notation
for positive sequences: $a_n \lesssim b_n$ if $\limsup_n a_n/b_n
<\infty$ and $a_n \asymp b_n$ if $a_n \lesssim b_n$ and $b_n \lesssim
a_n$ simultaneously.} $J_n$ with $2^{J_n} \asymp n$, the $\theta_{j,k}
$'s are drawn independently according to
%
\begin{equation}
\label{spikeslabinformal} \pi_j( dx ) = ( 1 - w_{j,n} )
\delta_0( dx) + w_{j,n} g (x)\,dx,
\end{equation}
for appropriate level-dependent weights $w_{j,n}>0$. Here, $\delta
_y(dx)$ is the Dirac mass at point $y$ and $g$ is a bounded density on
$\R$. For $j > J_n$, we put $\theta_{j,k}=0$. The construction of the
spike and slab prior does not involve knowledge of the smoothness index
$\beta$. Due to the point mass at zero, the posterior resembles many
properties of a wavelet thresholding procedure. In Theorem~\ref
{thWNlinf}, we prove adaptive concentration\vadjust{\goodbreak} rates
\[
\sup_{\theta_0 \in\Theta_\beta} E_{\theta_0}^n \bigl[
P^\pi \bigl(\theta: \llVert \theta-\theta_0\rrVert
_{L^\infty} \geq\epsilon_n(\beta)  \mid Y^n \bigr)
\bigr] \leq n^{-B}
\]
uniformly in $\beta\in\mathcal I$, where $\epsilon_n(\beta
)=M(n/\log n)^{-\beta/(2\beta+1)}$ and for some constants $B,M >0$
depending on $\pi$ and ${\mathcal I}$ only. Moreover, the polynomial
speed $n^{-B}$ at which the contraction holds is sharp according to
Theorem~\ref{thlowerbound} (up to the exponent~$B$).
The spike and slab prior (\ref{spikeslabinformal}) therefore leads
to an adaptive minimax posterior concentration rate over H\" older
balls $\mathcal H(\beta, L)$ for the sup-norm loss, without additional
$\log n$ term. To the best of our knowledge, this is the first
construction of a prior leading to an optimal adaptive posterior
concentration rate in sup-norm. However, we miss the optimal rate by a
logarithmic term if, for the same prior, we consider contraction under
the $L^2$-metric instead of $L^\infty$. We show in Theorem~\ref
{subsecWNL2} how to modify the spike and slab prior in order to
remove the logarithmic terms in the $L^2$-metric and achieve exact
adaptation in that setting too.

An answer to question (I) is presented in Section~\ref{secupper}. We
derive a generic upper bound, neither restricted to the white noise
model nor to $L^2$ or $L^\infty$ losses by considering priors which
are uniform over well chosen discrete sieves of $\Theta$.
In this abstract framework, Theorem~\ref{thuppbound} provides
conditions which imply that there exists a constant $C>0$ such that
%
\begin{equation}
\label{bornesupabstract} \sup_{\theta_0 \in\Theta} E_{\theta_0}^n
\bigl[P^\pi \bigl(\theta:\ell(\theta_0,\theta) >
\epsilon_n  \mid Y^n \bigr) \bigr] \leq e^{-Cn\Omega(\epsilon_n, \Theta, \ell)^2}.
\end{equation}
The interesting case is $n\Omega(\epsilon_n, \Theta, \ell
)^2\rightarrow\infty$, implying posterior concentration at rate
$\epsilon_n$. The rate can also be made adaptive by letting $\epsilon
_n = \epsilon_n(\theta_0)$ vary with $\theta_0$. Comparing (\ref
{bornesupabstract}) with the lower bound (\ref{lowerboundinformal}), we see in particular that the upper and lower bounds agree,
up to the constants $c$ and $C$, and are therefore sharp in that sense.
The rather abstract conditions which are required for (\ref{bornesupabstract}) are satisfied in the Gaussian white noise model and for
density estimation (Propositions~\ref{upbmodel}--\ref{teodensityestimation}).

In Section~\ref{discussion}, we discuss various implications of the
lower and upper bounds (\ref{lowerboundinformal}) and (\ref{bornesupabstract}). First, we outline how these bounds on posterior
concentration rates lead to the construction of Bayesian estimators
having asymptotic minimax (adaptive) frequentist risk, generalising the
result of Ghosal, Ghosh and van~der Vaart
\cite{ghosalghoshvdv00}, Theorem~2.5.
In Section~\ref{furtherremarks}, we point out the links between
posterior coverage and confidence balls. Interestingly, the lower bound
(\ref{lowerboundinformal}) implies that the classical strategy for
derivation of concentration rates fails if an arbitrary loss is
considered. This is developed in Section~\ref{subsecconseq}. Finally,
we discuss the relation of the derived results to other works in
Section~\ref{furtherbiblio}, both from a frequentist and Bayesian
point of view.

\section{A generic lower bound}\label{seclower} In this section, we
exhibit tractable conditions on the structure of a statistical experiment
$\mathcal E^n = \{P^n_\theta: \theta\in\Theta\}$ generated by data
$Y^n$ in order to obtain an explicit lower bound on
%
\[
E_{\theta_0}^n \bigl[ P^\pi\bigl(\theta:
\ell( \theta, \theta_0 ) > \epsilon_n \mid
Y^n\bigr) \bigr],\qquad\theta_0 \in\Theta,
\]
where $\epsilon_n$ can be either fixed or a function of $\theta_0$,
$\pi$ is a prior on $\Theta$, $P^\pi(\cdot\mid Y^n)$ denotes the
posterior distribution associated to $\pi$ and $\ell: \Theta\times
\Theta\rightarrow[0,\infty)$ is a given loss function.

Assume that
the parameter space $\Theta$ is equipped with a pre-metric $d$.
Let $\Theta_0 \subset\Theta$ and let $(\epsilon(\theta), \theta
\in\Theta_0)$ denote a collection of positive $\theta$-dependent
radii over $\Theta_0$.
We define a local and a global modulus of continuity related to
$\epsilon(\cdot)$ between $d$ and $\ell$ over a class $\Theta_0$ by setting
%
\begin{equation}
\label{modcontloc} \Omega \bigl(\epsilon(\cdot), \theta, \ell \bigr) = \inf\bigl\{ d
\bigl( \theta, \theta'\bigr):\ell\bigl(\theta, \theta'
\bigr) \geq\epsilon(\theta) + \epsilon\bigl(\theta'\bigr),
\theta' \in\Theta_0\bigr\}
\end{equation}
and
%
\begin{equation}
\label{modcontglob} \Omega\bigl(\epsilon(\cdot), \Theta_0, \ell\bigr) =
\inf_{\theta\in\Theta
_0} \Omega \bigl(\epsilon(\cdot), \theta, \ell \bigr).
\end{equation}
To illustrate the meaning of $\Omega$, consider, for instance, the
context of the Gaussian white noise model (\ref{WGN}) developed in
Section~\ref{secWN} below, where $\Theta_0 \subset\Theta=
L^2([0,1])$. Take $d= L^2$ and $\ell= L^\infty$ the sup-norm, and for
$\beta, L>0$, let $\Theta_0 = \mathcal H(\beta, L)$ be a H\"older
ball. Set $\epsilon_n(\theta)= M (n/\log n)^{-\beta/(2\beta+1)}$
for $\theta\in\Theta_0$ and $M>0$.
Then
%
\begin{equation}
\label{OmeganforHoelder} \Omega\bigl(\epsilon_n(\cdot), \theta,
L^\infty\bigr) \lesssim\sqrt{\log n /n}\qquad\mbox{for every }\theta\in
\Theta_0
\end{equation}
hence
\[
\Omega\bigl(\epsilon_n(\cdot), \Theta_0,
L^\infty\bigr)\lesssim\sqrt{\log n /n}
\]
as well (for a proof see Section~\ref{metricentropy}).
Similarly, if $\Theta_0 = \mathcal H(\beta_1, L) \supset\mathcal
H(\beta_2, L)$ with $\beta_1 <\beta_2$ and
if
%
\begin{equation}
\label{adaptiverate} \epsilon_n(\theta) = \cases{ M ( n/\log
n)^{-\beta_2 /(2\beta_2+1)}, &\quad if $\theta\in \mathcal H(\beta_2, L)$,
\cr
M
( n/\log n)^{-\beta_1 /(2\beta_1+1)}, &\quad otherwise,}
\end{equation}
then
%
\[
\Omega \bigl(\epsilon_n(\cdot), \Theta_0,
L^\infty \bigr) \lesssim \sqrt{\log n /n}.
\]
%
Obviously, when
$d= \ell$, then for all $\theta\in\Theta_0$ we have
$\Omega(\epsilon_n(\cdot), \theta, \ell) \geq\epsilon_n(\theta
)$ and $\Omega(\epsilon_n(\cdot), \Theta_0,\ell) \geq\inf\{
\epsilon_n(\theta), \theta\in\Theta_0\}$.

\begin{teo} \label{thlowerbound}
Let $\Theta_0 \subset\Theta$. Let $d$ be a pre-metric on $\Theta$.
Assume that $\ell$ is a pseudo-metric\footnote{That is, the
axioms of a metric are required with $\ell(\theta, \theta)=0$ but
possibly $\ell(\theta, \theta')=0$ for some distinct $\theta\neq
\theta'$.} on $\Theta_0$,
and that the prior $\pi$ and the family of positive sequences
$(\epsilon_n(\theta), \theta\in\Theta_0)$ satisfy the posterior
concentration condition:
%
\begin{equation}
\label{postrate} \sup_{\theta_0 \in\Theta_0} E_{\theta_0}^n
\bigl[P^\pi \bigl(\theta: \ell(\theta_0, \theta) \geq
\epsilon_n(\theta_0) \mid Y^n \bigr) \bigr]=
o(1).
\end{equation}
%
Assume that the family $\{P^n_\theta: \theta\in\Theta_0\}$ is
dominated by some $\sigma$-finite measure $\mu$ and that there exists
a constant $K>0$ such that
%
\begin{equation}
\label{A1} P_{\theta'}^n \bigl( \mathcal L_n
\bigl(\theta'\bigr) - \mathcal L_n(\theta) \geq K n d
\bigl(\theta, \theta'\bigr)^2 \bigr) = o(1),
\end{equation}
uniformly over all $\theta, \theta' \in\Theta_0$
satisfying
\[
\Omega\bigl(\epsilon_n(\cdot),\theta, \ell\bigr)\leq d\bigl(\theta,
\theta'\bigr) \leq2\Omega\bigl(\epsilon_n(\cdot),\theta,
\ell\bigr),
\]
where ${\mathcal L}_n(\theta) = \log\frac{dP^n_\theta}{d\mu
}(Y^n)$ denotes the log-likelihood function w.r.t. $\mu$. If $n \Omega
(\epsilon_n(\cdot),\Theta_0,\ell)^2\rightarrow\infty$, then, for
all $\theta_0 \in\Theta_0$ and large enough $n$
%
\begin{equation}
\label{lowerbound} E_{\theta_0}^n \bigl[ P^\pi \bigl(
\theta: \ell( \theta_0, \theta) > \epsilon_n(
\theta_0) \mid Y^n \bigr) \bigr] \geq e^{-3K n \Omega
(\epsilon_n(\cdot),\theta_0,\ell)^2}.
\end{equation}
\end{teo}

The proof is delayed until Section~\ref{secproofs}.

\begin{rem}
By taking $\epsilon_n(\theta) = \epsilon_n$ constant on $\Theta_0$,
we retrieve the more stringent result (\ref{lowerboundinformal})
announced in Section~\ref{mainresults}.
\end{rem}

%
\begin{rem}[(About the assumptions)]
Assumption (\ref{A1}) is merely on the pre-metric $d$ that must be
related to the intrinsic geometry of the experiment $\mathcal E^n$: it
shows in particular that $d$ must be able to control locally the
likelihood ratio. This can be the Hellinger distance used in the Birg\'
e--Le Cam testing theory in density estimation or simply the
$L^2$-distance in Gaussian white noise model linked to the Hilbert
space structure on which relies the existence of an iso-normal process.
Note also that since $d$ is not required to be symmetric, the order
$d(\theta, \theta')$ is important in assumption (\ref{A1}).
\end{rem}

In Sections~\ref{secWN} and~\ref{secupper}, we show that under some
additional assumptions the lower bound (\ref{lowerbound}) is sharp.

\section{Upper bounds in the white noise model via spike and slab priors} \label{secWN}
In this section, we prove that the lower bound obtained in (\ref
{lowerbound}) is sharp in the white noise model when $\ell$ is either
the sup-norm $L^\infty$ or the $L^2$-norm. This is done using spike
and slab type priors. We observe
%
\begin{equation}
\label{WGN} Y^n = \theta+ n^{-1/2}\dot W,
\end{equation}
where the signal of interest $\theta$ belongs to the Hilbert space
\[
\Theta= L^2\bigl([0,1]\bigr) = \biggl\{\theta:[0,1]\rightarrow\R
\mbox{ with } \int_{[0,1]}\theta(x)^2\,dx<\infty
\biggr\}
\]
and $\dot W$ is a Gaussian white noise on $\Theta$. The noise $\dot W$
is not realisable as a random element of $L^2$; it is therefore viewed
as the standard \textit{iso-Gaussian process} for the Hilbert space
$\Theta$.
Picking an orthonormal wavelet basis, we equivalently observe
%
\begin{eqnarray}
\label{WNwave} Y_{j,k} & =& \theta_{j,k} + n^{-1/2}
\epsilon_{j,k}, \qquad\epsilon _{j,k} \sim_{\mathrm{i.i.d.}}
\mathcal N(0,1), \qquad j \in\N,   k \in I_j,
\end{eqnarray}
where $\theta_{j,k}=\int_0^1\theta(x)\Psi_{j,k}(x)\,dx$ is the
wavelet coefficient associated to a given compactly supported wavelet
basis $(\Psi_{j,k})_{(j,k)\in\Lambda}$ of $\Theta$ with $\Lambda=
\{(j,k), k \in I_j, j\in\N\}$.
We append the basis with boundary conditions and assume that it is
associated with a $R$-regular multi-resolution of $L^2([0,1])$; see
\cite{cohen2000} and \cite{cdv1998}. The terms corresponding to $j=0$
incorporate the scaling function and we have that $\llvert  I_j\rrvert  $ is of order
$2^j$. We identify $\Theta= L^2([0,1])$ with
\[
\ell^2(\Lambda) = \biggl\{\theta= (\theta_{j,k})_{(j,k)\in\Lambda
}:
\sum_{(j,k)\in\Lambda}\theta_{j,k}^2<\infty
\biggr\}
\]
and we transfer two loss functions on the sequence space model:
the $L^2$-loss
\[
\ell_2\bigl(\theta, \theta'\bigr) = \biggl(\sum
_{(j,k)\in\Lambda}\bigl(\theta _{j,k}-
\theta_{j,k}'\bigr)^2 \biggr)^{1/2},
\]
and the $L^\infty$-loss
\[
\ell_\infty\bigl(\theta,\theta'\bigr)=\sum
_{j \in\N}2^{j/2}\max_{k\in
I_j}\bigl
\llvert \theta_{j,k}-\theta_{j,k}'\bigr\rrvert.
\]
%
Since $(\Psi_{j,k})_{(j,k)\in\Lambda}$ is orthonormal, $\ell_2$
coincides with the $L^2([0,1])$ norm. However, the losses $\ell_\infty
$ and $L^\infty$ are not comparable on $\Theta= L^2([0,1])$
identified with $\ell^2(\Lambda)$, but rather on smooth subspaces of
$\Theta$. To that end, introduce the H\"older balls\footnote{Having
$\beta= m+\{\beta\}$ with $m$ and integer and $\beta\in(0,1]$, the
class $\mathcal H(\beta, L)$ coincides with functions $f = \sum_{(j,k)\in\Lambda} \theta_{j,k}\psi_{j,k}$ that are $m$-times
differentiable with $f^{(m)}$ being H\"older continuous of order $\{
\beta\}$ provided the regularity of the multi-resolution exceeds
$\beta$.}
%
\begin{equation}
\label{eqdefHoelder} \mathcal H(\beta, L)= \bigl\{ \theta= (\theta_{j,k})_{(j,k)\in
\Lambda}
: \llvert \theta_{jk} \rrvert \leq L2^{-j(\beta+1 /2)}, (j,k)\in\Lambda
\bigr\}
\end{equation}
for $\beta>0, L>0$.
Then we also have that $\ell_{\infty}(\theta,\theta')$ and $\llVert
\theta-\theta'\rrVert  _{L^\infty([0,1])}$ are comparable on $\mathcal
H(\beta, L) \subset\ell^2(\Lambda)$.

\subsection{Adaptive posterior concentration rates under sup-norm loss: Spike and slab prior} \label{subsecWNsup}
Throughout the following, let $g$ be a bounded density on $\R$, which satisfies
\[
\inf_{x \in[-L_0, L_0]}g(x)>0
\]
for some $L_0 >0$.
We consider the following family of priors on $\Theta= \ell^2(\Lambda
)$. Set $J_n = \lfloor\log n/\log2 \rfloor$ and notice that $n/2 <
2^{J_n} \leq n$. For $j \leq J_n $ and $k \in I_j$, the $\theta_{j,k}
$'s are drawn independently from
%
\begin{eqnarray}
\label{spikeslab} \pi_j( dx ) &=& ( 1 - w_{j,n} )
\delta_0( dx) + w_{j,n} g (x)\,dx.
\end{eqnarray}
For $j>J_n$, $\pi_j(dx)= \delta_0(dx)$, or equivalently, $\theta
_{j,k}=0$. We assume that there are constants $K> 0$, $\tau>1/2$, such
that $n^{-K} \leq w_{j,n} \leq2^{-j(1+ \tau)}$, for all $j\leq J_n$.
This constraint on the mixture weights implies in particular that the
prior favours sparse models since the individual probability to be nonnull becomes small as the resolution level $j$ increases. We then have
the following.

\begin{teo} \label{thWNlinf}
Consider a prior in the family of spike and slab priors defined above.
If $Y^n$ is drawn from the white noise model (\ref{WNwave}), for
every $0<\beta_1 \leq\beta_2$ and $L_0-1 \geq L>0$, there exist $M,
B>0$ such that
\[
\sup_{\theta_0 \in\mathcal H(\beta, L) } E_{\theta_0}^n \bigl[
P^\pi \bigl( \theta: \ell_\infty( \theta,
\theta_0) \geq M (n/\log n)^{-\beta/ (2\beta+1)} \mid Y^n \bigr)
\bigr] \leq n^{-B}
\]
uniformly in $\beta\in[\beta_1, \beta_2]$.
\end{teo}

The proof of Theorem~\ref{thWNlinf} is given in Section~\ref
{prthWNlinf}. It is based on a fine description of the asymptotic
behaviour of the posterior distribution on the selected sets of
coefficients $\theta_{j,k}$, of the form $S = \{ (j,k), \theta_{j,k}
\neq0\}$, that is, we consider coefficients that are not equal
to $0$ under the posterior distribution. Lemma~\ref{lemWNSc} in
Section~\ref{prthWNlinf} states that the posterior distribution is
asymptotically neither forgetting nonnegligible coefficients $\theta
_{j,k}^{0}$ nor selecting too small coefficients $\theta_{j,k}^{0}$
under $P_{\theta_0}^n$ with $\theta_0=(\theta_{j,k}^0)_{(j,k)\in
\Lambda}$. As follows from the proof, if the prior density $g$ is
positive and continuous on $\R$, then the conclusion of Theorem~\ref
{thWNlinf} remains valid for every $L >0$ and the procedure is
independent of both the smoothness $\beta$ and the radius~$L$.

\begin{rem} Setting $\epsilon_n(\beta) = M (n/\log n)^{-\beta
/(2\beta+1)}$, we have
\[
\Omega\bigl(\epsilon_n(\beta), \mathcal H(\beta, L),
\ell_\infty\bigr) = O(\sqrt{\log n/n})
\]
and according to Theorem~\ref{thlowerbound}, the best possible
expectation of the posterior probability of complements on $\epsilon
_n(\beta)$ neighbourhoods in $\ell_\infty$ is at most of polynomial
order $n^{-B'}$ for some $B'>0$. Thus, Theorem~\ref{thWNlinf} is
sharp up to the constants $B'$ and~$B$.
\end{rem}

\subsection{Adaptive posterior concentration rates under $L^2$ loss: Block spike and slab prior} \label{subsecWNL2}
Theorem~\ref{thWNlinf} implies the existence of $\widetilde M>0$
such that
%
\begin{equation}
\label{WNL2local} E_{\theta_0}^n \bigl[ P^\pi \bigl(
\theta: \ell_2(\theta, \theta_0) \geq\widetilde M (n/\log
n)^{-\beta/ (2\beta+1)} \mid Y^n \bigr) \bigr] \leq n^{-B}
\end{equation}
uniformly in $\beta\in[\beta_1, \beta_2]$ since $\ell_2$ is
dominated by $\ell_\infty$. Therefore, an adaptive minimax posterior
concentration rate in the $\ell_2$-norm is also obtained up to a $\log
n$ term.
It can indeed be proved that for this prior the $\log n$ term cannot be
avoided. Since the spike and slab prior (\ref{spikeslab}) is a
product measure on the wavelet coefficients, this might be viewed as a
Bayesian analogue of the fact that separable rules do not give
adaptation with the clean rates in $\ell_2$ (cf. Cai \cite
{Cai2008421}). To circumvent this drawback, we propose a block spike
and slab prior which achieves the minimax adaptive rate for the $\ell
_2$-loss without additional $\log n$ term. The posterior associated to
this prior is easier to simulate from numerical data since the space of
possible selected sets is much smaller than the local spike and slab
prior (\ref{spikeslab}). It leads, however, to suboptimal posterior
concentration rates under sup-norm loss.

For $j \leq J_n$, pick a family of independent random vectors
$\underline{\theta}_j = (\theta_{j,k})_{k \in I_j}$ for $j \in\N$
according to the distribution
%
\begin{eqnarray}
\label{spikeslabL2} \widetilde\pi_j( dx ) &=& ( 1+ \nu_{j,n}
)^{-1}\bigl( \delta_0( dx) + \nu _{j,n}
g_j (x) \,dx\bigr) \qquad\forall x\in\R^{\llvert  I_j\rrvert},
\end{eqnarray}
where $g_j$ is a density with respect to Lebesgue measure on $\R
^{\llvert  I_j\rrvert  }$ which satisfies
%
\begin{equation}
\label{condgj} \sup_{x \in\R^{\llvert  I_j\rrvert  }}g_j(x) \leq
e^{ G \llvert  I_j\rrvert  }, \qquad\inf_{x\in
[-L_0,L_0]^{\llvert  I_j\rrvert  }} g_j(x) \geq
e^{-G\llvert  I_j\rrvert  },
\end{equation}
and $\nu_{j,n} = n^{\llvert  I_j\rrvert  /2} e^{-c\llvert  I_j\rrvert  }$ for some constants $G>0$ and
$c\geq4+G$. For $j>J_n$ put $\theta_{j,k}=0$. Condition (\ref
{condgj}) is satisfied in particular if, given that group $j$ is not
0, the $\theta_{j,k}$'s are i.i.d. with density $g$ satisfying the
same conditions as in the local spike and slab prior (\ref{spikeslab}).

\begin{teo} \label{thWNL2}
Consider a prior in the family of spike and slab priors defined above.
If $Y^n$ is drawn from the white noise model (\ref{WNwave}), for
every $0<\beta_1 \leq\beta_2$ and $L_0-1 \geq L>0$, there exist $M,
B>0$ such that
%
\[
\sup_{\theta_0 \in\mathcal H(\beta, L)} E_{\theta_0}^n \bigl[
P^\pi \bigl( \theta: \ell_2(\theta,\theta_0)
\geq M n^{-\beta/ (2\beta+1)} \mid Y^n \bigr) \bigr] \leq e^{-B
n^{1/(2\beta+1)}}
\]
uniformly in $\beta\in[\beta_1, \beta_2]$.
\end{teo}

The proof is given in Section~\ref{prthWNL2}.

\begin{rem}
Note that not only do we recover the optimal posterior concentration
rate (without any $\log n$ term) but we also bound from above the
expectation of the posterior concentration rate by a term of the order
\[
\exp \bigl(-cn\Omega\bigl(\zeta_n(\beta), \mathcal H(\beta, L), \ell
_2\bigr)^2 \bigr)
\]
with
$\zeta_n(\beta)= n^{-\beta/(2\beta+1)}$ when $\Omega$ is computed
under the intrinsic metric $d = \ell_2$. The same rate is provided by
the lower bound in Theorem~\ref{thlowerbound} and is therefore sharp
up to the constant $c>0$.
\end{rem}

\begin{rem}
Since the prior (\ref{spikeslabL2}) depends neither on $\beta$ nor
on $L$ (in particular if $g_j $ corresponds to $\llvert  I_j\rrvert  $ identically
distributed random variables with positive and continuous density $g$
on $\R$) the posterior concentration rate obtained in Theorem~\ref
{thWNL2} is moreover adaptive in the minimax sense of (\ref{defpostadap}).
\end{rem}

\section{A generic upper bound} \label{secupper}

In this section, we explore a more general situation and show that the
generic lower bound obtained in Theorem~\ref{thlowerbound} is indeed
sharp in a wider sense than the one considered in Section~\ref
{secWN}. In the context of an arbitrary experiment $\mathcal E^n$,
we construct priors with finite and increasing support, usually
referred to as sieve priors. Sieve priors have already been considered
in the Bayesian nonparametric literature in some specific context; see
\cite{ghosalghoshvdv00} and \cite{ghosalvaart2006}. In both cases,
the interest of these priors is that they lead to optimal posterior
concentration rates, without additional $\log n$ terms.
From a practical point of view, however, the construction of their
support and their implementation is close to being impossible.
Moreover, they have poor behaviour in terms of credible and confidence
sets. In this section, we shall use such priors in the same way, as a
device for the existence of an optimal estimation procedure, not as
priors to be used in practice.

We adopt the same framework as in Section~\ref{seclower}: $\mathcal
E^n = \{P_\theta^n,\theta\in\Theta\}$ is generated by the
observation $Y^n$ and is dominated by some $\sigma$-finite measure
$\mu$, and $\mathcal L_n(\theta) = \frac{dP_\theta^n}{d\mu}(Y^n)$
is a likelihood function. The loss function $\ell: \Theta\times
\Theta\rightarrow[0,\infty)$ is a pseudo-metric. Let us be given a family
\[
\epsilon_n= \bigl(\epsilon_n(\theta), \theta\in\Theta
\bigr)
\]
that we understand as the target posterior concentration rate at point
$\theta$ relative to the loss $\ell$.
Typically, $\epsilon_n(\theta)$ is the minimax rate of estimation
over a subclass $\Theta_0 \subset\Theta$ which contains $\theta$.
Let $(\Theta_n)_{n \geq1}$ be an increasing sequence of compact
subsets of $\Theta$ for the topology induced by the loss $\ell$. More
precisely, we only require that $\Theta_n$ can be covered by a finite
collection of balls centered at $\theta$ with radius $\epsilon
_n(\theta)$ in terms of the loss $\ell$.
We denote by $N_n$ the number of such balls and by $ \theta_{(l)}$ the
centers of these balls for $l=1,\ldots, N_n$. Note that we do not
necessarily require that $N_n$ is the minimal number of such balls
satisfying the coverage property. We define a \textit{sieve prior} as follows:
%
\begin{equation}
\label{sieveprior} \pi_n = \frac{1}{N_n} \sum
_{l=1}^{N_n} \delta_{ \theta_{(l)}}.
\end{equation}
To control the posterior concentration rate, we need to partition the
sieve $( \theta_{(l)}, 1 \leq l\leq N_n)$ into slices.

\begin{defi}
For every $\theta_0 \in\Theta_n$, a partition $(\mathcal J_r,1 \leq
r \leq R_n)$ of $\{ 1,\break  \ldots, N_n\}$ (we omit the dependence upon
$\theta_0$ in the notation) is called $\theta_0$-admissible if:
\begin{longlist}[(ii)]
\item[(i)] There exists $A>0$ such that $\mathcal J_0 = \{ l: \ell(
\theta_0, \theta_{(l)}) \leq A \epsilon_n(\theta_0)\}$.
\item[(ii)] For all $1 \leq r \leq R_n$, $\llvert  \mathcal J_r\rrvert   \leq
\llvert  \mathcal J_0\rrvert   $.
\end{longlist}
\end{defi}

\begin{teo}\label{thuppbound}
Assume that there exist constants $C_0, K_0, K_1>0$ and for every
$\theta_0\in\Theta_n$ a $\theta_0$-admissible partition $(\mathcal
J_r, 1 \leq r \leq R_n)$ together with injective maps $j_r: \mathcal
J_r\rightarrow\mathcal J_0$ such that
%
\begin{eqnarray}
\label{sievecond1} 
&&  P_{\theta_0}^n
\bigl(\exists r, \exists l: \mathcal L_n( \theta _{(l)}) -
\mathcal L_n(\theta_{(j_r(l))}) > -K_0 n d(
\theta_{(l)}, \theta_{(j_r(l))})^2 \bigr)
\nonumber\\[-8pt]\\[-8pt]\nonumber
&&\qquad \leq  C_0 \exp \bigl(-K_1 n \Omega\bigl(
\epsilon_n(\cdot), \theta _0,\ell\bigr)^2
\bigr)
\end{eqnarray}
and
%
\begin{equation}
\label{sievecond2} \sum_{r=1}^{R_n}e^{ -K_0 n u_r^2}
\leq C_0 e^{- K_1 n \Omega(\epsilon
_n(\cdot),\theta_0,\ell)^2},
\end{equation}
where $u_r = \min\{ d( \theta_{(l)}, \theta_{(j_r(l))}), l \in
\mathcal J_r\}$.
Then for all $\theta_0 \in\Theta_n$,
%
\begin{equation}
\label{upbound} E_{\theta_0}^n \bigl[ P^{\pi_n} \bigl(
\theta: \ell( \theta, \theta_0) > A \epsilon_n(
\theta_0)\mid Y^n \bigr) \bigr] \leq 2C_0
e^{ - K_1n \Omega(\epsilon_n(\cdot), \theta_0, \ell)^2}.
\end{equation}
\end{teo}

The proof is delayed until Section~\ref{secproofs}. Conditions (\ref
{sievecond1}) and (\ref{sievecond2}) on the admissible partition are
rather abstract. Interestingly, (\ref{sievecond1}) is the only
condition which links the geometry of $\Theta$ to the model $\{
P_\theta^n,\theta\in\Theta\}$. To illustrate conditions (\ref
{sievecond1})--(\ref{sievecond2}) and the admissible partition,
consider the following setup:
\[
\Theta= \bigcup_{\beta\in[\beta_1, \beta_2]} \mathcal H(\beta, L)=\mathcal
H(\beta_1, L) \subset\ell_2(\Lambda),
\]
where the H\"older ball $\mathcal H(\beta,L)$ is defined in
Section~\ref{secWN}.
Put $\ell(\theta, \theta') = \ell_\infty(\theta, \theta')$ and
$d(\theta, \theta') = \ell_2(\theta, \theta')$.
Let $\Theta_n = \{ \theta\in\Theta: \theta_{j,k}= 0,  \forall j >
J_n\}$ with $n < 2^{J_n}\leq2n$ and
set $\phi_n = \phi_0(\log n/n)^{1/2}$, where $\phi_0>0$ is fixed. Define
%
\begin{equation}
\label{Dnsievedef}\mathcal D_n = \bigl\{ \theta= (a_{j,k}
\phi_n, j\leq J_n, k\in I_j),
a_{j,k} \in\Z\cap[-L-1, L+1]\bigr\}
\end{equation}
and identify $\mathcal D_n$ as a subset of $\Theta_n$ by appending
zeros, that is, $\theta_{j,k}=0$ whenever $j>J_n$. The set ${\mathcal
D}_n$ defines the sieve, which we can enumerate as $\{ \theta_{(l)}, 1
\leq l \leq N_n\}$ with $N_n = \llvert  \mathcal D_n\rrvert  $. 
For any $\theta_0 = (\theta_{j,k}^0)_{(j,k) \in\Lambda} \in
\mathcal H(\beta, L)$ with $\beta\in[\beta_1, \beta_2]$, there
exists an integer $J_n(\beta)$ and a constant $b_0$ such that
%
\begin{eqnarray}
\label{eqJndefproperties}
\sup_{j >J_n(\beta)}\max_{k \in I_j} \bigl
\llvert \theta_{j,k}^0\bigr\rrvert &\leq&\phi
_n/4, \qquad\sum_{j>J_n(\beta)}2^{j/2}
\max_{k\in I_j} \bigl\llvert \theta _{j,k}^0
\bigr\rrvert \leq\epsilon_n(\beta),
\nonumber\\[-8pt]\\[-8pt]\nonumber
2^{J_n(\beta)} &\leq& b_0 (n/\log n)^{1/(2\beta+1)},
\end{eqnarray}
%
and we can pick $\theta^* \in\mathcal D_n$ satisfying
\begin{eqnarray*}
&&\forall(j,k)\in\Lambda, \qquad\bigl\llvert \theta_{j,k}^0
- \theta^*_{j,k}\bigr\rrvert \leq\phi_n/2.
\end{eqnarray*}
This\vspace*{1pt} implies in particular that $\forall j> J_n(\beta), \forall k \in
I_j, \theta_{j,k}^*=0$, and $\ell_\infty(\theta_0,\theta^*) \leq
(\phi_0+2) \epsilon_n(\beta)$.
We are ready to construct an admissible partition. First, consider the
semi-metric $d_1: \mathcal D_n\times\mathcal D_n \rightarrow
[0,\infty)$ defined by
\[
d_1\bigl(\theta, \theta'\bigr)^2 = \sum
_{j\leq J_n, k \in I_j} \bigl(\bigl(\theta _{j,k} -
\theta'_{j,k}\bigr)^2 - \phi_n^2
{\mathbf1}_{\{(j,k)\in
{\mathcal U}^c, \theta_{j,k} \wedge\theta'_{j,k}< \theta
_{j,k}^0<\theta_{j,k} \vee\theta'_{j,k}\}} \bigr),
\]
where
%
\begin{equation}
\label{defUinconstrofsieve} \mathcal U = \Bigl\{ (j,k) \in\Lambda, j\leq J_n,
\min_{t \in\Z
} \bigl\llvert \theta_{j,k}^0 -
t \phi_n\bigr\rrvert \leq\phi_n/4 \Bigr\}.
\end{equation}
Using the semi-metric $d_1$, we say that $\theta, \theta' \in
\mathcal D_n$ are equivalent if $d_1(\theta, \theta') = 0$, which
defines an equivalence relation. Denote by $\mathcal I_r$ the elements
of the corresponding quotient space, and let $\mathcal I_0$ be the
equivalence class of $\theta^*$. Then, for any $\theta\in\mathcal I_0$
\[
\ell_\infty(\theta_0, \theta) \leq \frac{3\phi_n}{4}\sum
_{j=0}^{J_n(\beta)} 2^{j/2} +
\epsilon_n(\beta) \leq\bigl(3\phi _0b_0^{1/2}
+1\bigr)\epsilon_n(\beta).
\]
Eventually, we can define for $A\geq4(3\phi_0b_0^{1/2}+1)$ the sets
%
\begin{equation}
\label{defpartition} \mathcal J_0 = \bigl\{ \theta\in\mathcal
D_n: \ell_\infty(\theta, \theta_0) \leq A
\epsilon_n(\beta)\bigr\}, \qquad\mathcal J_r = \mathcal
I_r \cap\mathcal J_0^c,
\end{equation}
where we have identified the partition of the indices with the
partition of the elements of $\mathcal D_n$.
We then have the following. 

\begin{prop} \label{upbmodel}
Assume that $\theta_0 \in\bigcup_{\beta\in[\beta_1,\beta_2
]}\mathcal H(\beta,L)=\mathcal H(\beta_1,L)$, and consider the
partition $(\mathcal J_r, r\geq0)$ (depending on $\theta_0$) defined
as in (\ref{defpartition}) above. Then, if $\ell=\ell_\infty$, the
partition $(\mathcal J_r, r\geq0)$ is $\theta_0$-admissible and
satisfies (\ref{sievecond2}).

Moreover, if $Y^n$ is drawn from the white noise model (\ref
{WNwave}), for every $0<\beta_1 \leq\beta_2$ and $L>0$, there exist
$M, B>0$ such that
\[
\sup_{\theta_0 \in\mathcal H(\beta, L) } E_{\theta_0}^n \bigl[
P^{\pi_n} \bigl( \theta: \ell_\infty( \theta,
\theta_0) \geq M (n/\log n)^{-\beta/ (2\beta+1)} \mid Y^n \bigr)
\bigr] \leq n^{-B}
\]
uniformly in $\beta\in[\beta_1, \beta_2]$.
\end{prop}

The proof of Proposition~\ref{upbmodel} is given in Appendix~\ref{appprlemupb}.
The generic upper bound allows us to prove posterior\vspace*{1pt} concentration in
$L^2$ loss with the ``clean'' adaptive rate $\epsilon_n(\beta
)=n^{-\beta/(2\beta+1)}$ as well. In fact, we obtain an analogous
result to Theorem~\ref{thWNL2} by constructing an appropriate sieve
prior and using Theorem~\ref{thuppbound}. For sake of brevity, we
give the statement without a proof.

\begin{prop} \label{updiscreteWNL2}
There exists a sieve prior $\pi_n$, such that if $Y^n$ is drawn from
the white noise model (\ref{WNwave}), for every $0<\beta_1 \leq
\beta_2$ and $L>0$, there exist $M, B>0$ with
\[
\sup_{\theta_0 \in\mathcal H(\beta, L) } E_{\theta_0}^n \bigl[
P^{\pi_n} \bigl( \theta: \ell_2 ( \theta,
\theta_0) \geq M n^{-\beta/ (2\beta+1)} \mid Y^n \bigr) \bigr]
\leq\exp \bigl(-Bn^{1/(2\beta+1)}\bigr)
\]
uniformly in $\beta\in[\beta_1, \beta_2]$.
\end{prop}

Even more interesting is that the generic upper bound can be also
applied to prove adaptive rates for density estimation, with respect to
$\ell_\infty$ loss. In this model, we observe $Y^n = (Y_1,\ldots, Y_n)$,
where $Y_i, i=1,\ldots,n$ are independent and identically distributed
on $[0,1]$ with density $f_\theta$ and write
%
\begin{equation}
\label{eqsqrtdensexpansion} \sqrt{f_\theta(x)} = \sum_{(j,k)\in\Lambda}
\theta_{j,k} \Psi_{j,k}(x).
\end{equation}
Here, the parameter space consists of vectors $\theta= (\theta
_{j,k})_{(j,k)\in\Lambda} \in\mathcal H(\beta, L)$ such that the
right-hand side of (\ref{eqsqrtdensexpansion}) is larger than some
constant $c>0$ and $\llVert   \theta\rrVert  _{\ell_2}=1$. We refer to this
restricted H\"older space as $\mathcal H'(\beta,L)$ in the sequel. In
this case, we can take $d=\ell_2$ again.

\begin{prop}\label{teodensityestimation}
There exists a sieve prior $\pi_n$, such that if $Y^n$ is drawn from
the density model (\ref{eqsqrtdensexpansion}), for every $1/2<\beta
_1 \leq\beta_2$ and $L>0$, there exist $M, B>0$ with
\[
\sup_{\theta_0 \in\mathcal H'(\beta, L) } E_{\theta_0}^n \bigl[
P^{\pi_n} \bigl( \theta: \ell_\infty( \theta,
\theta_0) \geq M (n/\log n)^{-\beta/ (2\beta+1)} \mid Y^n \bigr)
\bigr] \leq n^{-B}
\]
uniformly in $\beta\in[\beta_1, \beta_2]$.
\end{prop}

The proof of Proposition~\ref{teodensityestimation} is given in
Section~\ref{appthmthmupbdensityestimation}.

\section{Further results and discussion} \label{discussion}

\subsection{Construction of minimax adaptive estimators given adaptive concentration} \label{freq-consequences}

The main focus of this work is to study the full posterior distribution
under the frequentist assumption of a true parameter $\theta_0$. As a
statistical implication of the results let us shortly comment on
convergence rates of Bayesian point estimators. In the nonadaptive
case, Theorem~2.5 in \cite{ghosalghoshvdv00} asserts the existence of
an estimator that converges with the posterior concentration rate to
the true parameter. However, the construction of the estimator
crucially depends on knowledge of the rate $\epsilon_n$ and is
therefore not applicable in the adaptive setup. Not surprisingly, the
Bayes estimator
%
\begin{equation}
\label{bayesestimator} \widehat\theta\in\mathop{\operatorname{argmin}}_{\delta}
E^{\pi} \bigl[\ell (\delta,\theta) \mid Y^n \bigr],
\end{equation}
(assumed to be well defined) (see, e.g., \cite{robert2004}), Chapter~2,
will achieve the adaptive rate under quite general conditions. To see
this, assume that\break $\ell(\theta, \theta') = \ell(\theta',\theta)$
for all $\theta, \theta' \in\Theta$ and observe that for any
$\theta_0\in\Theta$,
%
\begin{eqnarray}
\ell(\widehat\theta, \theta_0) &\leq& E^{\pi}\bigl[\ell(
\widehat\theta, \theta)+\ell(\theta,\theta_0)\mid Y^n
\bigr] \leq2E^{\pi}\bigl[\ell(\theta,\theta_0)\mid
Y^n\bigr]. \label{eqriskineq}
\end{eqnarray}
If the loss is bounded, say $\sup_{\theta\in\Theta}\ell(\theta,
\theta_0)\leq M$, this can be further controlled by
\begin{eqnarray*}
&& 2 \bigl( \epsilon_n(\beta) +M P^\pi \bigl(\theta: \ell(
\theta, \theta_0)>\epsilon_n(\beta) \mid Y^n
\bigr) \bigr).
\end{eqnarray*}
Consider now a subset $\Theta_\beta\subset\Theta$ such that for any
$\theta_0 \in\Theta_\beta$ the posterior concentration rate at
$\theta_0$ is bounded by $\epsilon_n(\beta)$ in the slightly
stricter sense
\[
\sup_{\theta_0\in\Theta_\beta} E_{\theta_0}^n \bigl[
P^\pi \bigl(\theta: \ell(\theta, \theta_0)>
\epsilon_n(\beta) \mid Y^n \bigr) \bigr]= o\bigl(
\epsilon_n(\beta) \bigr).
\]
Then
\[
\sup_{\theta_0\in\Theta_\beta} P_{\theta_0}^n \bigl(
\ell (\widehat\theta, \theta_0) > 2(M+1) \epsilon_n(\beta)
\bigr) = o(1)
\]
and
\[
\sup_{\theta_0\in\Theta_\beta} E_{\theta_0}^n \bigl[
\ell (\widehat\theta, \theta_0) \bigr] = O\bigl(\epsilon_n(
\beta) \bigr).
\]
Consequently, $\widehat\theta$ achieves the rate $\epsilon_n(\beta
)$ over $\Theta_\beta$. In the case of an unbounded loss functions
$\ell$, slightly refined arguments can be applied. Consider, for
instance, the framework of Theorem~\ref{thWNL2}. Here, adaptation is
meant over H\"older balls $\mathcal H(\beta, L)\subset\mathcal
H(\beta_1, L)$ with $\beta_1>0$. Since $\sup_{\theta,\theta'\in
\mathcal H(\beta_1, L)}\ell_2(\theta, \theta')\leq M_2 <\infty$
for some constant $M_2$, we can improve any estimator by projection on
$\mathcal H(\beta_1, L)$. The projected estimator lies then within
$\ell_2$-distance $M_2$ from $\theta_0$. Thus, considering risk of
estimators, we may replace the $\ell_2$-loss by the modified bounded
loss function $\widetilde\ell_2 = \min(\ell_2, M_2)$. Together with
Theorem~\ref{thWNL2} and the steps described above, the Bayes
estimator with respect to $\widetilde\ell_2$ yields then an adaptive
estimator.

An alternative modification to incorporate unbounded loss functions
goes via a slicing of $\ell(\theta, \theta_0)$ in $E_{\theta_0}^n
[E^{\pi}[\ell(\theta,\theta_0)\mid Y^n]]$. With (\ref{eqriskineq}),
\begin{eqnarray*}
&& E_{\theta_0}^n \bigl[\ell(\widehat\theta, \theta_0)
\bigr] \leq2\epsilon_n(\beta) +2\sum_{j\geq1}
(j+1)\epsilon_n(\beta) E_{\theta_0}^n \bigl[
P^\pi\bigl(\theta: \ell(\theta,\theta_0)> j
\epsilon_n(\beta) \mid Y^n \bigr) \bigr].
\end{eqnarray*}
The second term of the upper bound will typically be negligible
(uniformly over~$\Theta_\beta$) since it involves the posterior concentration.
In fact, under the conditions of Theorem~\ref{thWNlinf}, the Bayes
estimator in (\ref{bayesestimator}) adapts to H\"older balls with
respect to the sup-norm loss.

\begin{prop}\label{proprisksup}
Consider the spike and slab prior (\ref{spikeslab}) with $w_{j,n}\leq
n^{-6}2^{-j(1+\tau)}$ and $\tau>1/2$. If $Y^n$ is drawn from the
white noise model (\ref{WNwave}), for any $0<\beta_1 \leq\beta_2$
and\vspace*{1pt} $L_0-1 \geq L>0$, then
there exists $M>0$ such that with $\epsilon_n(\beta)=M (n/\log
n)^{-\beta/ (2\beta+1)}$,
\[
\sup_{\beta\in[ \beta_1, \beta_2]}  \sup_{\theta_0 \in\mathcal H(\beta, L) }
P_{\theta_0}^n \bigl( \ell _\infty(\widehat\theta,
\theta_0) \geq \epsilon_n(\beta) \bigr) = o(1)
\]
and
\[
\sup_{\beta\in[ \beta_1, \beta_2]} \epsilon_n(\beta)^{-1} \sup
_{\theta_0 \in\mathcal H(\beta, L) } E_{\theta_0}^n \bigl[ \ell
_\infty(\widehat\theta, \theta_0) \bigr] <\infty.
\]
\end{prop}

The proof of Proposition~\ref{proprisksup} is given together with
the proof of Theorem~\ref{thWNlinf} in Section~\ref{prthWNlinf}.

\subsection{Posterior concentration and confidence balls} \label{furtherremarks}
The posterior distribution does not only provide point estimators but
also Bayesian measures of uncertainty. Apart from regular parametric
models, it is not clear whether such credible sets have a frequentist
interpretation as measures of confidence. In this section we discuss
some consequences on the asymptotic behaviour of posterior credible balls.

Assume that the prior $\pi$ leads to a concentration rate $\epsilon
_n$ over some subset $\Theta_0$ of the parameter space, that is,
\[
\sup_{\theta_0\in\Theta_0}   E_{\theta_0}^n \bigl[
P^\pi \bigl(\theta: \ell( \theta, \theta_0) >
\epsilon_n \mid Y^n \bigr) \bigr] \leq e^{ - c n \Omega(\epsilon_n, \Theta_0, \ell
)^2}
\rightarrow0. %
\]
As discussed in Section~\ref{freq-consequences}, this implies under
mild conditions existence of a point estimator $\widehat\theta_n$ satisfying
\[
\sup_{\theta_0\in\Theta_0}   P_{\theta_0}^n \bigl( \ell(
\widehat\theta _n, \theta_0) > \epsilon_n \bigr) =
o(1).
\]
Let $\alpha_n\in(0,1)$ be a sequence, possibly tending to zero, that
satisfies\break $e^{ - c n \Omega(\epsilon_n, \Theta_0, \ell)^2}
= o(\alpha_n)$. Construct the credible ball
\[
C_n = \bigl\{ \theta: \ell( \theta, \widehat\theta_n) \leq
q_{\alpha
_n}^\pi\bigr\}, %
\]
where $q_{\alpha_n}^\pi$ is the $1 - \alpha_n$ posterior quantile of
$\ell(\theta, \widehat\theta_n)$ so that
%
\begin{equation}
\label{baycoverage} P^\pi \bigl( \theta\in C_n \mid
Y^n \bigr) \geq1 - \alpha_n.
\end{equation}
We then have the following two properties for $C_n$:
%
\begin{eqnarray}
\label{conf1} \int_\Theta P_\theta^n(
\theta\in C_n )\,d\pi(\theta) &\geq& 1 - \alpha_n,
\nonumber\\[-8pt]\\[-8pt]\nonumber
\sup_{\theta\in\Theta_0} P_\theta^n\bigl(
\ell(C_n) > 4\epsilon_n \bigr) &=& o(1),
\end{eqnarray}
where $\ell(C_n) = \sup\{ \ell(\theta, \theta'): \theta, \theta
' \in C_n \} = 2 q_{\alpha_n}^\pi$ is the diameter of $C_n$.

\begin{pf*}{Proof of (\ref{conf1})}
The first inequality is a consequence of the Fubini theorem since (\ref
{baycoverage}) is true for all $Y^n$ so that
\[
\int_\Theta P_\theta^n( \theta\in
C_n )\,d\pi(\theta) = \int_{\mathcal Y^n } P^\pi
\bigl( \theta\in C_n \mid Y^n \bigr)\,d m_\pi
\bigl(Y^n\bigr) \geq1 - \alpha_n,
\]
where $m_\pi$ is the marginal distribution of $Y^n$.
The second statement of (\ref{conf1}) follows from the fact that for
all $\theta\in C_n$,
\[
\ell( \widehat\theta_n, \theta) \geq-\ell( \theta_0, \widehat
\theta _n) + \ell( \theta_0, \theta).
\]
Thus, on the event $\{\ell( \theta_0, \widehat\theta_n)\leq\epsilon
_n\}$, for all $t < q_{\alpha_n}^\pi$ and every $\theta_0 \in\Theta_0$,
\[
\alpha_n \leq P^\pi\bigl(\ell( \theta, \widehat
\theta_n)> t \mid Y^n \bigr) \leq P^\pi \bigl(
\ell( \theta_0, \theta)> t - \epsilon_n \mid
Y^n \bigr)
\]
implying
\begin{eqnarray*}
P_{\theta_0}^n\bigl( q_{\alpha_n}^\pi> 2
\epsilon_n \bigr) &\leq& P_{\theta_0}^n\bigl(\ell(
\theta_0, \widehat\theta_n)> \epsilon_n\bigr)+
P_{\theta_0}^n\bigl( P^\pi \bigl( \ell(
\theta_0, \theta)> \epsilon_n \mid Y^n \bigr)
\geq \alpha_n \bigr)
\\
&=& o(1) + \frac{ e^{-c n \Omega(\epsilon_n, \Theta_0, \ell)^2}}{
\alpha_n} = o(1),
\end{eqnarray*}
uniformly over $\theta_0\in\Theta_0$. This completes the proof of
(\ref{conf1}).
\end{pf*}

A natural question is then whether the first inequality of (\ref
{conf1}) can be turned into
\[
\inf_{\theta\in\Theta} P_\theta^n( \theta\in
C_n ) \geq1 - \alpha
\]
at least for some reasonably small $\alpha$. Of particular interest is
the case of adaptive posterior concentration rate, which we illustrate
considering the sup-norm loss $\ell_\infty$ over a collection of H\"
older balls $\bigcup_{\beta\in[\beta_1,\beta_2]} \mathcal H(\beta, L
) = \mathcal H(\beta_1, L)$ with $0<\beta_1 \leq\beta_2$ and $L>0$ fixed.
Assume that
\[
\sup_{\beta\in[ \beta_1, \beta_2]}  \sup_{\theta_0 \in\mathcal
H(\beta,L) }
E_{\theta_0}^n \bigl[ P^\pi \bigl(\theta:
\ell_\infty( \theta, \theta_0) > \epsilon_n(
\beta) \mid Y^n \bigr) \bigr] \leq n^{-B} %
\]
with $\epsilon_n(\beta) = M (n/\log n)^{-\beta/(2\beta+1)}$ for
some $M, B>0$ and
\[
\sup_{\beta\in[ \beta_1, \beta_2]}  \sup_{\theta_0 \in\mathcal H(\beta, L) }
P_{\theta_0}^n \bigl( \ell _\infty(\widehat\theta,
\theta_0) \geq \epsilon_n(\beta) \bigr) = o(1).
\]
By Theorem~\ref{thWNlinf} and Proposition~\ref{proprisksup}, this
is, for instance, achieved by the prior in~(\ref{spikeslab}) and the
Bayes estimator $\widehat\theta$. Let $\alpha_n \geq n^{-B+ t}$ for
some $t>0$, then following from (\ref{conf1}) we obtain
%
\begin{eqnarray}
\label{conf2}  \int_\Theta P_\theta^n(
\theta\in C_n )\,d\pi(\theta) &\geq& 1 - \alpha_n,
\nonumber\\[-8pt]\\[-8pt]\nonumber
 \sup_{\beta\in[\beta_1,\beta_2]}  \sup_{\theta_0 \in\mathcal
H(\beta, L) }
P_{\theta_0}^n\bigl( \ell(C_n) > 2
\epsilon_n(\beta) \bigr) &=&o(1).
\end{eqnarray}
In this case, there exists no adaptive confidence band
(see, e.g., \cite{HoffmannNickl11}), so that (\ref{conf2}) implies that
for all $\alpha>0$,
\[
\lim_n \inf_{\beta\in[\beta_1, \beta_2]}  \inf
_{\theta_0 \in
\mathcal H(\beta, L) } P_{\theta_0}^n( \theta_0
\in C_n) = 0.
\]
The nonexistence of adaptive confidence bands means that requiring
both honest frequentist coverage and adaptive length of the band is too
strong. Integrating out the confidence band is a weaker notion and a
possible alternative to the approach of \cite{HoffmannNickl11} which
modifies confidence bands by taking off some points $\theta$ and by
demanding coverage and adaptive length over this restricted set.
Further notice that the first inequality of (\ref{conf2}) implies that
\[
\pi \bigl( \theta: P_\theta^n( \theta\in C_n )
\leq1 - \alpha \bigr) \leq\frac{ \alpha_n }{ \alpha}.
\]
This is, however, not enough to characterise the parameter values
$\theta_0$ for which
$P_{\theta_0}^n( \theta_0 \in C_n)$ is small. This question is of
interest but beyond the scope of the present paper.

\subsection{Consequences on proving strategies for posterior concentration rates} \label{subsecconseq}
The lower bound in Theorem~\ref{thlowerbound} has an interesting
consequence for nonparametric Bayes in general. So far, the state of
the art techniques for deriving posterior consistency and concentration
rates date back to the work of \cite{schwartz}. Her approach relies on
two key ideas. First, treat the numerator and denominator in the Bayes
formula separately. Second, introduce an abstract test and express the
upper bound in terms of errors of the first and second type. These
methods were later refined by \cite{barron1988,ghosalghoshvdv00} and \cite{ghosalvaart2006}. In particular, from the
proof of Theorem 1 of \cite{ghosalvaart2006}, if for $\epsilon_n$
their conditions (2.4), (2.5) (associated to the loss $\ell$) are
satisfied and
\[
P_{\theta_0}^n \bigl[ \mathcal L_n(\theta) -
\mathcal L_n(\theta_0) < -n \epsilon_n^2
\bigr] \leq e^{-c_1 n\epsilon_n^2}
\]
for some positive $c_1$, then
\[
E_{\theta_0}^n \bigl[ P^\pi \bigl( \theta: \ell(
\theta, \theta_0)> M\epsilon_n \mid Y^n \bigr)
\bigr] \lesssim e^{-
c_2 n \epsilon_n^2 }
\]
for some $c_2>0$. The lower bound of Theorem~\ref{thlowerbound},
however, implies that
\[
\Omega( \epsilon_n, \theta_0, \ell) \gtrsim
\epsilon_n.
\]
Therefore, if the targeted concentration rate (say the minimax
estimation rate over some given class) $\epsilon_n^*$ satisfies
\[
\Omega\bigl(\epsilon_n^*, \theta_0, \ell\bigr) = o\bigl(
\epsilon_n^*\bigr)
\]
then the approach of \cite{ghosalvaart2006}, Theorem 1, leads to a
suboptimal posterior concentration rate.
The core of the problem comes from the decomposition of the posterior
probability which treats separately the denominator $D_n$ and the
numerator $N_n$ (see Section~\ref{subsecderivationofinequalofggvdv}) where the main steps of the arguments of \cite
{ghosalvaart2006} are recalled. Denote by $\Phi_{n}$ the test for
$H_0: \theta=\theta_0$ versus $H_1: \ell( \theta, \theta
_0)>\epsilon_n$. Then the derived upper bound can be written as
follows: There exists a positive sequence $u_n$ such that
%
\begin{eqnarray}
\label{equpperboundheuristic}
&& E_{\theta_0}^n \bigl[ P^\pi \bigl(
\theta: \ell(\theta_0,\theta )>\epsilon_n \mid
Y^n \bigr) \bigr]
\nonumber\\[-8pt]\\[-8pt]\nonumber
&& \qquad \leq   E_{\theta_0}^n [\Phi_{n} ] +
e^{cnu_n^2}\sup_{\theta: \ell(\theta_0,\theta)>\epsilon
_n}E_{\theta}^n [1-
\Phi_{n} ]+e^{-c'nu_n^2},
\end{eqnarray}
with finite constants $c,c'>0$ on which we do not have good control.
For the right-hand term of (\ref{equpperboundheuristic}) to be
small, we need
\[
\sup_{\theta\in\Theta_n} {\mathbf1}_{\{\ell( \theta_0, \theta)
>\epsilon_n\}} E_\theta^n[
1 - \Phi_n]= o\bigl(e^{-cnu_n^2} \bigr).
\]
Hence, $\epsilon_n$ shall verify the constraint $\Omega( \epsilon_n,
\theta_0, \ell) \gtrsim u_n $; if the minimax estimation rate
$\epsilon_n^*$ over a given class satisfies
$\Omega( \epsilon_n^*, \theta_0, \ell) = o(u_n) $, the approach
through tests typically leads to suboptimal posterior concentration rates.
To illustrate this, consider the white noise model where $d$ is the
$L^2$ loss, $\ell= \ell_\infty$, and $\theta_0$ belonging to a H\"
older ball with smoothness $\beta$. Assume that $\theta_1 \in\Theta
$ satisfies $\ell(\theta_0,\theta_1)>\epsilon_n$ and
$\llVert  \theta_1 - \theta_0 \rrVert  _{L^2} \leq C \Omega( \epsilon_n(\cdot),
\theta_0, \ell) $ for some fixed arbitrary $C$. Any test $\Phi_n$
with error of first kind smaller than some small $\epsilon$ must have
a second kind error greater than that of the likelihood ratio test
$\phi_{n,\theta_1}$ for $H_0:\theta=\theta_0$ against $H_1:\theta
=\theta_1$. In other words,
\[
E_{\theta_1}^n[ 1 - \Phi_n] \geq
E_{\theta_1}^n[ 1 - \phi_{n,\theta
_1}] \gtrsim
e^{- c_1 n \llVert  \theta_1 -\theta_0\rrVert  _{L^2}^2} \geq e^{ -
nc_1 C \Omega( \epsilon_n, \theta_0, \ell)^2 }
\]
for some $c_1 >0$. This implies $ \Omega( \epsilon_n, \theta_0, \ell
) \gtrsim u_n$. The above argument can be generalised to other models,
in particular, to density estimation with $\theta_0 \in\mathcal
H(\beta, L)$ for $L, \beta>0$. If we rely on the bound (\ref{equpperboundheuristic}), the best achievable concentration rate is
given by the $(\ell, d)$-modulus of continuity $\omega(u_n)$ as
defined in (\ref{modcontdef}) below.
As an example, consider density estimation. Any\vspace*{1pt} prior which leads to
the minimax estimation error $n^{-\beta/(2\beta+1)}$ in the Hellinger
metric gives $u_n = n^{-\beta/(2 \beta+ 1)}$ in (\ref{equpperboundheuristic}) (possibly\vspace*{1pt} up to $\log n$ terms). Since for $\ell$ the
sup-norm and $\theta_0\in\mathcal H(\beta,L)$, $\omega(n^{-\beta
/(2\beta+1)}) \lesssim n^{-(\beta-1/2)/(2 \beta+1)}$, this explains
the (suboptimal) rate observed in \cite{GineNickl11} which was derived
using the standard approach, and thus a bound of the type~(\ref{equpperboundheuristic}).

\subsection{Relation to other works} \label{furtherbiblio}
In the last decade, a variety of posterior concentration rates have
been derived. These studies include density estimation in the case of
independent and identically distributed observations as in \cite
{ghosalghoshvdv00}, nonparametric regression (Ghosal and van~der Vaart \cite{ghosalvaart2006})
and the white noise model (Zhao \cite{zhao00}, Belitser and
Ghosal \cite{belitserghosal03}), Markov
models (Tang and Ghosal \cite{tangghosal07}), Gaussian time series (Choudhuri, Ghosal and
Roy \cite{chgh04} and Rousseau, Chopin and Liseo
\cite{rcl}) to name but a few, or the recent canonical statistical
setting of \cite{castillokp2013}. For each of these models, a variety
of families of priors have been investigated.
An interesting feature of the Bayesian nonparametric approaches
considered in these papers is that minimax adaptive concentration rates
are achieved using hierarchical types of priors where, at the highest
level of hierarchy some hyperparameter, somehow related to the class
index $\beta$, is itself given a prior distribution. For instance, in
the case of density estimation, the renown class of Dirichlet process
mixtures or related types of mixtures lead to adaptive posterior
concentration rates over collections of H\"older balls of regularity
$\beta$, up to a $\log n $ term, see, for instance, \cite{rousseau09,kruijerrousseauvdv10,ghosalshentokdar} and \cite
{scricciolo12} under the Hellinger or the $L^1$ losses on the densities.
Gaussian random fields, with inverse Gamma bandwidth as prior models
also lead to adaptive posterior concentration rates up to a $\log n$
term for a large class of models, including the nonlinear regression
model under the empirical quadratic loss on the design and the
classification problem under the $L^2$ loss; see \cite{vzantenvdv09}.
Similarly, orthonormal basis expansions with random truncation
generically yield adaptive posterior concentration rates too, provided
the loss function is well chosen; see \cite{arbeletal13}. All these
results, however, are proved using the approach proposed by \cite
{ghosalvaart2006}, which relies on the existence of tests with
exponentially small error of the second kind outside $\ell
$-neighbourhoods of the true parameter. Therefore, these results are
applicable to loss functions which behave similarly to $d$.

Previous to this work, suboptimal asymptotic behaviour of posterior
distributions has been observed for specific loss functions.
Arbel, Gayraud and and Rousseau
\cite{arbeletal13} shows, for instance, that a random truncation prior with
minimax adaptive (up to a $\log n$ term) posterior concentration rate
under $L^2$ loss leads to significantly suboptimal posterior
concentration rate and (and risk) under pointwise loss.

To our knowledge, the question of the existence of adaptive minimax
posterior concentration rates when $\ell$ is the pointwise loss or
even the sup-norm loss $L^\infty$ has been an open question until now.
A consequence of our results is the explicit construction of priors
that lead to adaptive concentration rates for various loss functions
(including the sup-norm $L^\infty$). Given a prior $\pi$ and a loss
function $\ell$, the best achievable rate of concentration of the
posterior distribution is intimately linked to the geometry of the
experiment $\mathcal E^n = \{P_\theta^n,\theta\in\Theta\}$ in the
most classical sense of Le Cam (see, e.g., \cite
{LeCamYang86}), expressed through the pre-metric $d$. The behaviour of
such a pair $(\ell, d)$ is reminiscent of several phenomena in minimax
theory: these include estimation of linear functionals \cite
{donoholiu}, constrained risk inequalities \cite
{brownlow96,CaiLow05,CaiLowZhao04} or the existence of adaptive
confidence sets \cite{Low97,cailow04}. In all these studies, a key
ingredient is the behaviour of a $(\ell, d)$ modulus of continuity
%
\begin{equation}
\label{modcontdef} \omega(\epsilon) = \sup\bigl\{\ell\bigl(\theta,
\theta'\bigr): d\bigl(\theta, \theta'\bigr)\leq
\epsilon, \theta, \theta' \in\Theta\bigr\},\qquad \epsilon>0
\end{equation}
that quantifies the maximal error in the desired $\ell$-loss for a
prescribed statistical distance $\epsilon$ induced by the experiment
$\mathcal E^n$ via the intrinsic pre-metric $d$. More precisely, if
there are two sequences $\epsilon_n>0$ and $\theta_n \in\Theta$,
such that $d(\theta_0, \theta_n) \leq\epsilon_n$ implies that there
exists no convergent test of
\[
H_0:\theta= \theta_0\quad\mbox{against}\quad
H_n:\theta= \theta_n,
\]
then $\omega(\epsilon_n)$ yields a lower bound for the minimax
estimation rate of $\theta$ in $\ell$-loss.
The nonexistence of adaptive confidence intervals over H\"older balls
in the Gaussian white noise model lies at the heart of this simple
phenomenon: In that case, $\ell$ is the pointwise or $L^\infty$-norm
and $d$ is the $L^2$-metric. The fact that an irregular function can be
close to a smooth functions in $L^2([0,1])$ while away from the smooth
target in $L^\infty$ explains the negative result of Low \cite{Low97}
[see also \cite{GineNickl11,HoffmannNickl11}] and is
quantified by $\omega(\epsilon_n)$.
Interestingly, in the Bayesian framework, the lower bound derived in
Theorem~\ref{thlowerbound} is of similar nature and the modulus of
continuity defined in (\ref{modcontglob}) is the dual of the modulus
of continuity considered in the frequentist minimax literature and
defined in (\ref{modcontdef}).

\section{Proofs} \label{secproofs}

\subsection{Proof of Theorem \texorpdfstring{\protect\ref{thlowerbound}}{2.1}}
We prove Theorem~\ref{thlowerbound} by contradiction. Assume that
there exist $\theta_0\in\Theta_0$ such that
%
\begin{equation}
\label{nolb} E_{\theta_0}^n \bigl[ P^\pi \bigl(
\theta: \ell(\theta_0, \theta) > \epsilon_n(
\theta_0) \mid Y^n \bigr) \bigr] < e^{- 3K n \Omega
(\epsilon_n(\cdot),\theta_0,\ell)^2},
\end{equation}
infinitely often, which without loss of generality we can assume to be
satisfied for all $n$. By definition of $\Omega(\epsilon_n(\cdot
),\theta_0,\ell)$, we can choose a sequence
$(\theta_n^*)_n \subset\Theta_0$ satisfying
\[
\Omega\bigl(\epsilon_n(\cdot), \theta_0, \ell\bigr) \leq
d\bigl(\theta_0, \theta _n^*\bigr) \leq2\Omega\bigl(
\epsilon_n(\cdot),\theta_0,\ell\bigr)
\]
and
\[
\ell\bigl( \theta_0, \theta_n^*\bigr) \geq
\epsilon_n(\theta_0) + \epsilon _n\bigl(
\theta_n^*\bigr)
\]
simultaneously. Then for every $\theta\in\Theta$,
\[
\ell\bigl( \theta, \theta_n^*\bigr) < \epsilon_n\bigl(
\theta_n^*\bigr)\quad \Rightarrow\quad \ell( \theta, \theta_0) >
\ell\bigl(\theta_0, \theta_n^*\bigr) - \epsilon
_n\bigl(\theta_n^*\bigr) \geq\epsilon_n(
\theta_0)
\]
so that
\begin{eqnarray*}
&&  E_{\theta_n^*}^n \bigl[ P^\pi \bigl(\theta: \ell
\bigl(\theta_n^*, \theta\bigr)< \epsilon_n\bigl(
\theta_n^*\bigr) \mid Y^n \bigr) \bigr]
\\
&&\qquad \leq  E_{\theta_n^*}^n \bigl[ P^\pi \bigl(\theta:
\ell( \theta _0, \theta)> \epsilon_n(\theta_0)
\mid Y^n \bigr) \bigr]
\\
&&\qquad \leq   e^{K n d(\theta_0, \theta_n^*)^2} E_{\theta_0}^n \bigl[
P^\pi \bigl(\theta: \ell( \theta_0, \theta)>
\epsilon_n(\theta _0) \mid Y^n \bigr) \bigr]
\\
&&\quad\qquad{} + P_{\theta_n^*}^n \bigl(\mathcal L_n\bigl(
\theta_n^*\bigr) - \mathcal L_n(\theta_0) > Kn
d\bigl(\theta_0, \theta_n^*\bigr)^2 \bigr)
\\
&&\qquad \leq  e^{-(3K - 2K) n \Omega(\epsilon_n(\cdot), \theta_0,\ell
)^2 } + P_{\theta_n^*}^n \bigl(\mathcal
L_n\bigl(\theta_n^*\bigr) - \mathcal L_n(
\theta_0) > Kn d\bigl(\theta_0, \theta_n^*
\bigr)^2 \bigr)
\\
&&\qquad =    o(1)
\end{eqnarray*}
in contradiction with the posterior concentration (\ref{postrate}).

\subsection{Proof of Theorem \texorpdfstring{\protect\ref{thuppbound}}{4.1}}
For $\theta_0 \in\Theta_n$, let $(\mathcal J_r, 0\leq r \leq R_n)$
be a $\theta_0$-ad\-missible partition satisfying (\ref{sievecond1})
and (\ref{sievecond2}). Let $A_n(\theta_0) = \{ \theta\in\Theta
_n: \ell(\theta_0, \theta) > A \epsilon_n(\theta_0)\} $, where $A$
is defined via the admissible partition in (\ref{defpartition}).
Set\break $p_{n,\theta}(Y^n)=\frac{dP_\theta^n}{d\mu}(Y^n)$ so that
$\mathcal L_n(\theta) = \log p_{n,\theta}$.
For the sieve prior $\pi_n$ defined in (\ref{sieveprior}),
\begin{eqnarray*}
P^{\pi_n} \bigl( A_n(\theta_0) \mid
Y^n \bigr) &=&\frac{ \sum_{l=1}^{N_n} {\mathbf
1}_{A_n(\theta_0)} ( \theta_{(l)})p_{n, \theta_{(l)}}(Y^n) }{\sum_{l=1}^{N_n} p_{n, \theta_{(l)}}(Y^n) }
\\
&\leq&\sum_{r=1}^{R_n} \frac{ \sum_{l \in\mathcal J_r} p_{n, \theta
_{(l)}}(Y^n) }{\sum_{l \in\mathcal J_0} p_{n, \theta_{(l)}}(Y^n) }
\leq \sum_{r=1}^{R_n} \max
_{l \in\mathcal J_r} e^{ \mathcal L_n(
\theta_{(l)}) -\mathcal L_n( \theta_{(j_r(l) )}) }.
\end{eqnarray*}
Let
\[
\Omega_n(\theta_0) = \bigl\{ \forall r \geq1, \forall l
\in \mathcal J_r: \mathcal L_n( \theta_{(l)}) -
\mathcal L_n(\theta _{(j_r(l))}) \leq-K_0 n d(
\theta_{(l)}, \theta_{(j_r(l))})^2 \bigr\}.
\]
On $\{Y^n \in\Omega_n(\theta_0)\}$,
\begin{eqnarray*}
P^{\pi_n} \bigl( A_n(\theta_0) \mid
Y^n \bigr) &\leq& \sum_{r=1}^{R_n}
e^{ - K_0 n d( \theta_{(l)}, \theta_{(j_r(l))})^2}
\\
&\leq& \sum_{r=1}^{R_n} e^{-K_0n u_r^2}
\leq C_0 e^{ - K_1 n \Omega
(\epsilon_n(\cdot), \theta_0, \ell)^2 },
\end{eqnarray*}
thanks to assumption (\ref{sievecond2}), which combined with
assumption (\ref{sievecond1}) completes the proof.
\subsection{Proof of Theorem \texorpdfstring{\protect\ref{thWNlinf}}{3.1} and of Proposition \texorpdfstring{\protect\ref{proprisksup}}{4}} \label{prthWNlinf}

Recall that the prior can be written in the following hierarchical way:
First, select a set of nonzero components $S$, with distribution $P$
the product of independent Bernoulli random variables $\mathcal B(
w_{j,n})$ for $j \leq J_n$. Given $S$, draw independently
$\theta_{j,k} \sim g$ for all $(j,k) \in S$, and put $\theta_{j,k}=0$
otherwise.

Asymptotically, the posterior concentrates on supports $S$ containing
only indices $(j,k)$ with $\llvert  \theta_{j,k}^0\rrvert  >\underline\gamma\sqrt
{\log n/n}$ and all indices $(j,k)$ with $\llvert  \theta_{j,k}^0\rrvert  >\overline
\gamma\sqrt{\log n/n}$, where $0<\underline\gamma<\overline\gamma
<\infty$ are appropriate constants. In this respect, the posterior
behaves similar as hard thresholding. Indeed, for
\[
\mathcal J_n(\gamma) = \bigl\{ (j,k)\in\Lambda: \bigl\llvert \theta
_{j,k}^0\bigr\rrvert >\gamma\sqrt{\log n/n} \bigr\}\qquad
\mbox{with }\gamma>0
\]
we have the following.

\begin{lemma}\label{lemWNSc}
Under the conditions of Theorem~\ref{thWNlinf}, for every $0<\beta
_1 \leq\beta_2$, $L\leq L_0-1$, and any $B>0$, there exists $\overline
\gamma>0$ such that
%
\begin{equation}
\label{WNSc1} \sup_{\beta_1 \leq\beta\leq\beta_2}  \sup_{\theta_0 \in\mathcal H(\beta,L)}
E_{\theta_0}^n \bigl[ P^\pi \bigl(S^c
\cap\mathcal J_n(\overline\gamma)\neq\varnothing\mid Y^n
\bigr) \bigr] \lesssim\frac{\log n}{n^B}.
\end{equation}
Suppose that the mixing weights in the spike and slab prior (\ref
{spikeslab}) satisfy $w_{j,n}\leq\min(\frac{1}2, n^{(\tau\wedge
1)/2-1/4-2B}2^{-j(1+\tau)})$ with $B> 0, \tau>1/2$. Then, for
sufficiently small $0<\underline\gamma$,
%
\begin{equation}
\label{WNSc2} \sup_{\beta_1 \leq\beta\leq\beta_2}  \sup_{\theta_0 \in\mathcal H(\beta, L) }
E_{\theta_0}^n \bigl[ P^\pi \bigl( S \cap\mathcal
J_n^c(\underline\gamma) \neq \varnothing\mid
Y^n \bigr) \bigr] \lesssim\frac{\log n}{n^B}.
\end{equation}
\end{lemma}

The proof of Lemma~\ref{lemWNSc} is delayed until Appendix~\ref
{appWNSL}. Suppose for the moment that for any $B>0$ the bound
$\overline\gamma$ can be chosen large enough such that
%
\begin{equation}
\label{maximumunderposteriortoshow} \sup_{\theta_0 \in\mathcal H(\beta, L) } E_{\theta_0}^n
\Bigl[ P^\pi \Bigl( \max_{(j,k)\in\mathcal J_n(\underline\gamma)} \bigl\llvert \theta
_{j,k}-\theta_{j,k}^0\bigr\rrvert >\overline
\gamma\sqrt{\log n/n} \mid Y^n \Bigr) \Bigr] \lesssim
\frac{\log n}{n^B},
\end{equation}
uniformly in $\beta\in[\beta_1,\beta_2]$. The last estimate ensures
that the posterior concentrates around $\theta_{j,k}^0$ with\vspace*{1pt} the good
rate $\sqrt{\log n/n}$ on every component $(j,k)$ on which signal
might be detected. Now we are ready to complete the proof of
Theorem~\ref{thWNlinf}. The definition of a H\"older ball in (\ref{eqdefHoelder}) implies that there exists a $J_n(\beta)$ with $2^{J_n(\beta
)}\leq k_2(n/\log n)^{1/(2\beta+1)}$ for some constant $k_2>0$ such
that $\mathcal J_n(\underline\gamma)\subset\{(j,k): j\leq
J_n(\beta), k\in I_j\}$ and
\[
\sup_{\theta_0\in\mathcal H(\beta, L) } \sum_{j>J_n(\beta
)}2^{j/2}
\max_{k\in I_j} \bigl\llvert \theta_{j,k}^0
\bigr\rrvert \leq\frac{1}2 M (n/\log n)^{-\beta/(2\beta+1)} =:
\frac{1}2 \epsilon_n(\beta),
\]
for $M$ a sufficiently large constant. In order to prove the theorem,
it is sufficient to show that $\ell_\infty(\theta,\theta_0)\leq
\epsilon_n(\beta)$ for all $\theta$ with $\max_{(j,k)\in\mathcal
J_n(\underline\gamma)} \llvert  \theta_{j,k}-\theta_{j,k}^0\rrvert  \leq\overline
\gamma\sqrt{\log n/n}$ and support $S$ satisfying the constraints
$S^c\cap\mathcal J_n(\overline\gamma)=\varnothing$ and $S \cap
\mathcal J_n^c(\underline\gamma)= \varnothing$. Using the properties
of $J_n(\beta)$,
\[
\ell_\infty(\theta,\theta_0) \leq\sum
_{j=0}^{J_n(\beta)} 2^{j/2} \max
_{k\in I_j} \bigl\llvert \theta _{j,k}-
\theta_{j,k}^0\bigr\rrvert +\frac{1}2
\epsilon_n(\beta)\leq\overline \gamma2^{J_n(\beta)/2} \sqrt{\log n/n}+
\frac{1}2\epsilon_n(\beta)
\]
and the right-hand side can further be uniformly bounded by $\epsilon
_n(\beta)$. This establishes Theorem~\ref{thWNlinf} provided (\ref
{maximumunderposteriortoshow}) is true.

For Theorem~\ref{thWNlinf}, it therefore remains to show (\ref
{maximumunderposteriortoshow}). By Lemma~\ref{lemWNSc}, we can
restrict ourselves to parameters with support $S$ satisfying $S^c\cap
\mathcal J_n(\overline\gamma)= \varnothing$. Using a union bound and
considering the cases $\underline\gamma\sqrt{\log n/n}<\llvert  \theta
_{j,k}^0\rrvert  \leq\overline\gamma\sqrt{\log n/n}$ and $\overline\gamma
\sqrt{\log n/n}<\llvert  \theta_{j,k}^0\rrvert  $ separately,
\begin{eqnarray*}
&& P^\pi \Bigl( \max_{(j,k)\in\mathcal J_n(\underline\gamma)} \bigl\llvert
\theta_{j,k}-\theta_{j,k}^0\bigr\rrvert >\overline
\gamma\sqrt{\log n/n}\mbox{ and } S^c\cap\mathcal
J_n(\overline\gamma)= \varnothing \mid Y^n \Bigr)
\\
&&\qquad \lesssim  n \max_{(j,k)\in\mathcal J_n(\underline\gamma)} P^\pi \bigl( \bigl
\llvert \theta_{j,k}-\theta_{j,k}^0\bigr\rrvert >
\overline\gamma\sqrt{\log n/n}\mbox{ and }\theta_{j,k}\neq0 \mid
Y^n \bigr).
\end{eqnarray*}
Consider the event
\[
\Omega_{n,B} = \bigl\{\sqrt{n}\bigl\llvert Y_{j,k} -
\theta_{j,k}^0\bigr\rrvert \leq\bigl(2 \log \llvert
I_j\rrvert + 2B \log n\bigr)^{1/2}, \forall j \leq
J_{n}, \forall k \in I_j\bigr\}.
\]
Then
%
\begin{equation}
\label{Omegan} P_{\theta_0}^n \bigl( \Omega_{n,B}^c
\bigr) \leq2 n^{-B} J_{n} \leq \frac{2\log n}{n^B}.
\end{equation}
For all $j \leq J_{n}, k \in I_j$, on $\Omega_{n,B}$, $\llvert  Y_{j,k}\rrvert  \leq\llvert
\theta_{j,k}^0\rrvert   + \frac{1}{2}\leq L_0-\frac{1}2$ and so if $a = \inf
\{g(x):   \llvert  x\rrvert   \leq L_0\}>0$, then, setting $u_0 = \Phi^{-1}( (1+
1/\sqrt{2})/2)$
%
\begin{eqnarray}
\label{lowerbounddenominator}
\quad && \int_{\R} e^{-\sklfrac{n}{2}( Y_{j,k} - \theta)^2}g(\theta)\,d\theta
\geq a(2\pi/n)^{1/2}\bigl(2 \Phi(u_0)-1\bigr) \geq a(
\pi/n)^{1/2},
\end{eqnarray}
where $\Phi(x)=\operatorname{Pr}(\mathcal N(0,1)\leq x)$. For any $(j,k)
\in\mathcal J_n(\underline\gamma)$,
\begin{eqnarray*}
&& P^\pi \bigl( \bigl\llvert \theta_{j,k} -
\theta_{j,k}^0\bigr\rrvert >\overline\gamma \sqrt{\log n/n}
\mbox{ and }\theta_{j,k} \neq0  \mid Y^n \bigr)
\\
&&\qquad \leq  a^{-1}\sup_x g(x) \biggl(
\frac{n}{\pi} \biggr)^{1/2} \int_{\R
} {\mathbf1}
\bigl\{ \bigl\llvert \theta- \theta_{j,k}^0\bigr\rrvert >
\overline\gamma\sqrt{\log n/n} \bigr\} e^{- \sklfrac{n}{2} ( Y_{j,k} - \theta)^2} \,d \theta.
\end{eqnarray*}
On $\Omega_{n,B}$,
\[
\bigl\{ \bigl\llvert \theta- \theta_{j,k}^0\bigr\rrvert >
\overline\gamma\sqrt{\log n/n}\bigr\} \subset\bigl\{ \llvert \theta-
Y_{j,k}\rrvert >\tfrac{1}2 \overline \gamma\sqrt{\log n/n}\bigr\}
\]
provided $\overline\gamma$ is large enough. Therefore, for any $(j,k)
\in\mathcal J_n(\underline\gamma)$,
\[
P^\pi \bigl( \bigl\llvert \theta_{j,k} -
\theta_{j,k}^0\bigr\rrvert >\overline\gamma\sqrt {\log n/n}
\mbox{ and }\theta_{j,k} \neq0  \mid Y^n \bigr)
\leq a^{-1}\sup_x g(x) 2^{3/2}e^{-\overline\gamma^2 \log n/8}
\]
and together with the union bound and the\vspace*{1pt} estimate of $P_{\theta_0}^n(
\Omega_{n,B}^c )$ above, equation~(\ref{maximumunderposteriortoshow}) follows for $\overline\gamma$ sufficiently large. The proof of
Theorem~\ref{thWNlinf} is complete.

The proof of Proposition~\ref{proprisksup} relies on the
computations above. Define $ A_1 = \{\max_{(j,k)\in\mathcal
J_n(\underline\gamma)} \llvert  \theta_{j,k}-\theta_{j,k}^0\rrvert  \leq\overline
\gamma\sqrt{\log n/n} \}$, $A_2 =\{ S: S^c\cap\mathcal
J_n(\overline\gamma)= \varnothing\}$ and $A_3 =\{ S: S \cap\mathcal
J_n(\underline\gamma) = \varnothing\}$. On\vspace*{1pt} $A_1\cap A_2\cap A_3$,
$\ell_\infty( \theta, \theta_0) \leq M (n/\log n)^{-\beta/(2\beta
+1)}$ for some $M>0$. Thus, with (\ref{eqriskineq}), Proposition
\ref{proprisksup} is proved if
%
\begin{equation}
\label{maxineqrisk} E_{\theta_0} \bigl[ E^\pi \bigl(
\ell_\infty(\theta, \theta_0) ({\mathbf1}_{A_1^c \cap A_2\cap A_3}+{
\mathbf1}_{A_2^c} +{\mathbf1}_{A_3^c} ) \mid Y^n \bigr) \bigr]
\leq\frac{\log n}{\sqrt{n}}.
\end{equation}
Let $A $ be a measurable subset of the parameter set, then using the
Cauchy--Schwarz inequality twice,
\begin{eqnarray*}
&& E_{\theta_0}^n \bigl[ E^\pi \bigl(
\ell_\infty(\theta, \theta _0){\mathbf1}_A \mid
Y^n \bigr) \bigr]
\\
&&\qquad \lesssim \sum_{j,k}2^{j/2}
E_{\theta_0}^n \bigl[ E^\pi \bigl( \bigl\llvert
\theta_{j,k} - \theta_{j,k}^0\bigr\rrvert
^2 \mid Y^n \bigr) \bigr]^{1/2}
E_{\theta_0}^n \bigl[ P^\pi\bigl(A\mid
Y^n\bigr) \bigr]^{1/2}
\\
&&\qquad \leq2 \sum_{j,k}2^{j/2}
\biggl(E_{\theta_0}^n \bigl[ E^\pi \bigl( \llvert
\theta_{j,k} - Y_{j,k}\rrvert ^2 \mid
Y^n \bigr) \bigr]+\frac{1}n \biggr)^{1/2}
E_{\theta_0}^n \bigl[ P^\pi\bigl(A\mid
Y^n\bigr) \bigr]^{1/2}.
\end{eqnarray*}
We apply this inequality to $ A = A_1^c \cap A_2 \cap A_3$, $A_2^c$,
and $A_3^c$. Using the bounds above, it is\vspace*{1pt} sufficient to control
$E_{\theta_0}^n[ E^\pi( ( \theta_{j,k} - Y_{j,k})^2 \mid Y^n )]$.
Recall the definition of the spike and slab prior (\ref{spikeslab})
and observe
\begin{eqnarray*}
&& E^\pi \bigl( ( \theta_{j,k} - Y_{j,k})^2
\mid Y^n \bigr)
\\
&&\qquad \leq Y_{j,k}^2+\frac{2 w_{j,n} \sup_x g(x) }{a} \frac{ \int_{\R} (\theta-Y_{j,k})^2 e^{-n (\theta-Y_{j,k})^2/2}
\,d\theta}{ \int_{-L_0}^{L_0} e^{-n (\theta-Y_{j,k})^2/2} \,d\theta}
\\
&&\qquad \leq Y_{j,k}^2+ \frac{2 w_{j,n} \sup_x g(x) }{a n }
\\
&&\quad\qquad{}\times  \bigl[ \Phi\bigl(
\sqrt{n}\bigl(L_0 -\theta^0_{j,k}\bigr)-
\epsilon_{j,k} \bigr) - \Phi\bigl( - \sqrt {n}\bigl(L_0 +
\theta^0_{j,k}\bigr)- \epsilon_{j,k} \bigr)
\bigr]^{-1}
\end{eqnarray*}
with $\epsilon_{j,k} = \sqrt{n} (Y_{j,k}-\theta_{j,k}^0)$ and $\Phi
$ the distribution function of a standard normal random variable. Since
$\llvert  \theta_{j,k}^0\rrvert  \leq L_0-1$,
\begin{eqnarray*}
&& \int_0^{\infty} e^{-\epsilon^2/2} \bigl( \Phi\bigl(
\sqrt{n}\bigl(L_0 -\theta^0_{j,k}\bigr)-
\epsilon\bigr) - \Phi\bigl( - \sqrt{n}\bigl(L_0 +\theta
^0_{j,k}\bigr)- \epsilon\bigr) \bigr)^{-1}\,d
\epsilon
\\
&& \qquad \leq \int_0^{\infty} e^{-\epsilon^2/2} \bigl( \Phi(
\sqrt{n}- \epsilon) - \Phi( - \sqrt{n}- \epsilon) \bigr)^{-1}\,d\epsilon
\\
&&\qquad  \lesssim \int_0^{\sqrt{n}} e^{-\epsilon^2/2} \,d
\epsilon+ \int_{\sqrt{n}}^\infty\epsilon e^{n/2-\epsilon\sqrt{n} }\,d
\epsilon
\\
&&\qquad  \lesssim1 + e^{-n/2}.
\end{eqnarray*}
The same inequality can be obtained for the integral over $(-\infty,
0)$. Consequently, there exists a universal constant $C>0$ for which
\[
E_{\theta_0}^n \bigl[ E^\pi \bigl( (
\theta_{j,k} - Y_{j,k})^2 \mid Y^n
\bigr) \bigr]\leq\bigl(\theta_{j,k}^0\bigr)^2+
\frac{1}n+ \frac{2 Cw_{j,n}
\sup_x g(x) }{a}.
\]
Since $\llvert  \theta_{j,k}^0\rrvert  \lesssim 2^{-j/2}$ and $w_{j,n}\leq 2^{-j}$,
we obtain that for any measurable set $A$ and uniformly over $\theta_0
\in\mathcal H(\beta, L)$,
\[
E_{\theta_0}^n \bigl[ E^\pi \bigl(
\ell_\infty(\theta, \theta _0){\mathbf1}_A \mid
Y^n \bigr) \bigr] \lesssim n E_{\theta_0}^n \bigl[
P^\pi \bigl( A \mid Y^n \bigr) \bigr]^{1/2}.
\]
From the proof of Theorem~\ref{thWNlinf} above, we find that the
right-hand side is of order $\log n/\sqrt{n}$, provided that the
exponent $B$ in Lemma~\ref{lemWNSc} and (\ref{maximumunderposteriortoshow}) is three. This completes the proof of Proposition
\ref{proprisksup}.

\subsection{Proof of Theorem \texorpdfstring{\protect\ref{thWNL2}}{3.2}} \label{prthWNL2}
We set $\uY_j = (Y_{j,k}, k \in I_j)$ and similarly $\utheta_j =
(\theta_{j,k}, k \in I_j)$. Whenever convenient, we identify $\uY_j$
and $\utheta_j$ as sequences indexed by the whole set of indices
$\Lambda$, setting their value to be $0$ on the complement of $I_j$.
Thus, if $\llVert  \cdot\rrVert  $ denotes the usual Euclidean norm on $\R
^{\llvert  I_j\rrvert  }$, we have $\ell_2(\utheta_j, \utheta_j') = \llVert  \utheta
_j-\utheta_j'\rrVert  $ with a slight abuse of notation.\vspace*{1pt}

The proof of Theorem~\ref{thWNL2} follows the classical line for
studying posterior concentration rates as proposed in \cite
{ghosalvaart2006}, with some extra care that has to be taken in order
to avoid the usual $\log n$ term that appears in this case. Set
$u_n(\beta) = n^{-\beta/(2\beta+1)}$ and let $\widetilde J_n(\beta) $
satisfy $ K_1 n^{1/(2\beta+1)} \leq2^{\widetilde J_n(\beta)} \leq2K_1
n^{1/(2\beta+1)}$, where $K_1$ is a constant to be large enough.
Define
\[
\Theta_n(\beta) = \{ \theta: \theta_{j,k}=
\theta_{j,k} {\mathbf 1}_{\{j\leq\widetilde J_n(\beta), k \in I_j\}} \}.
\]
%
We first prove that for some $c_1, K_1>0$,
%
\begin{equation}
\label{Thetan} P^\pi \bigl( \Theta_n(\beta)^c
\mid Y^n \bigr) \leq e^{-c_1 n u_n^2(\beta)}.
\end{equation}
Let $\theta_0 = (\theta_{j,k}^0)_{(j,k)\in\Lambda} \in\mathcal
H(\beta,L)$ with $L\leq L_0-1$ and $L_0$ is the constant appearing in
condition (\ref{condgj}).
Denote by $B_n$ the intersection of the events
\[
\bigl\{ Y^n: n  \bigl\llVert \uY_j -
\utheta_j^0\bigr\rrVert ^2 \leq e \llvert
I_j\rrvert, \forall j \mbox{ with } \widetilde J_n(
\beta)\leq j \leq J_n\bigr\}
\]
and
\[
\bigl\{ Y^n: \bigl\llvert Y_{j,k}-\theta_{j,k}^0
\bigr\rrvert \leq1/2, \forall j\leq\widetilde J_n(\beta), k\in
I_j\bigr\}.
\]
Set $c_e = e/2-1$ and $C_e = (1- e^{-c_e})^{-1}$. For a $\chi_p^2$
distributed random variable $\xi$, we have $\operatorname{Pr}(\xi>eq)\leq
e^{-c_eq}$ whenever $q\geq p$. Hence,
\begin{eqnarray*}
&& P_{\theta_0}^n \bigl( B_n^c \bigr)
\leq2ne^{-n/8}+ \sum_{j =\widetilde
J_n(\beta)}^{J_n}
e^{-c_e \llvert  I_j\rrvert  } \leq2ne^{-n/8}+C_e e^{-c_e
\llvert  I_{\widetilde J_n(\beta)}\rrvert   } \leq
e^{ - An^{1/(2\beta+1)}}, 
\end{eqnarray*}
for $n$ large enough, with $A$ proportional to $K_1$. Since
\[
\Theta_n(\beta)^c = \bigcup
_{j \geq\widetilde J_n(\beta)} \{ \theta: \theta_{I_j} \neq0\}
\]
(here $\theta_{I_j} \neq0$ means $\theta_{j,k} \neq0$ for at least
one $k\in I_j$) we conclude
\begin{eqnarray*}
P^\pi \bigl( \Theta_n(\beta)^c \mid
Y^n \bigr) &\leq& \sum_{j \geq\widetilde J_n(\beta)}
\frac{(1+\nu_{j,n})^{-1} \nu_{j,n}
\int_{\R^{\llvert  I_j\rrvert  }} e^{-\sklfrac{ n}{2} \llVert  \utheta_j - \uY_j\rrVert  ^2
}g_j(\utheta_j) \,d\utheta_j }{ \int_{\R^{\llvert  I_j\rrvert  }} e^{-\sklfrac{ n}{2}
\llVert   \utheta_j - \uY_j\rrVert  ^2 }d\widetilde\pi_j(\utheta_j )}
\\
&\leq& \sum_{j \geq\widetilde J_n(\beta)}e^{G \llvert  I_j\rrvert  }
\nu_{j,n} \biggl( \frac{ 2 \pi}{n} \biggr)^{\llvert  I_j\rrvert  /2} \exp \biggl(
\frac{ n \llVert  \uY_j\rrVert
^2}{2} \biggr).
\end{eqnarray*}
%
For all $j \geq\widetilde J_n(\beta)$, we have $\llVert  \utheta_j^0\rrVert  ^2 \leq
L^2 \llvert  I_j\rrvert  2^{-j(2\beta+1)}\leq C\llvert  I_j\rrvert  /n $, for the radius of the H\"
older ball $L$ and some constant $C>0$ which decreases to zero as $K_1$
grows. On the event $B_n$, we thus infer $n\llVert   \uY_j\rrVert  ^2 \leq2(C+e)
\llvert  I_j\rrvert  $. Therefore, on $B_n$,
\begin{eqnarray*}
P^\pi \bigl( \Theta_n(\beta)^c \mid
Y^n \bigr) &\leq& \sum_{j \geq\widetilde J_n(\beta)}e^{(G + C+e) \llvert  I_j\rrvert  }
\nu_{j,n} \biggl( \frac{ 2 \pi}{n} \biggr)^{\llvert  I_j\rrvert  /2}
\\
&\leq& 2e^{ - (c-G- e-\sklfrac{1}2 \log2\pi-C) \llvert  I_{\widetilde J_n(\beta)}\rrvert  } \leq2e^{-b K_1 n^{1/(2\beta+1)} }
\end{eqnarray*}
for some $b>0$ as soon as $c> G+e+\frac{1}2 \log2\pi$ provided we
choose $K_1$ large enough. This proves (\ref{Thetan}).
We are ready to complete the proof. For $A_n=\{\theta: (\sum_{j=0}^{\widetilde J_n(\beta)} \llVert   \utheta_j - \utheta_j^0\rrVert  ^2)^{1/2}
\leq Mu_n(\beta)/2\}$,
\begin{eqnarray*}
&&  P^\pi \bigl( \theta: \ell_2(\theta,
\theta_0) > Mu_n(\beta) \mid Y^n \bigr)
\\
&&\qquad \leq  P^\pi \bigl(\bigl\{\theta: \ell_2(\theta,
\theta_0) > Mu_n(\beta) \bigr\} \cap\Theta_n(
\beta) \mid Y^n \bigr) + P^\pi \bigl(\Theta_n(
\beta)^c \mid Y^n\bigr)
\\
&&\qquad \leq  P^\pi\bigl( A_n^c\mid
Y^n \bigr) + P^\pi \bigl(\Theta_n(
\beta)^c \mid Y^n \bigr).
\end{eqnarray*}
We bound the first term of the right-hand side by
\begin{eqnarray*}
P^\pi\bigl( A_n^c\mid Y^n \bigr)
&\leq& \int_{A_n^c}\prod_{j=0}^{\widetilde J_n(\beta
)}\frac{ e^{-\sklfrac{ n}{2} \llVert   \utheta_j - \uY_j\rrVert  ^2 }(1+\nu
_{j,n})\,d\widetilde\pi_j(\utheta_j) }{ \int_{\R^{\llvert  I_j\rrvert  }} e^{-\sklfrac{ n}{2} \llVert   \utheta_j - \uY_j\rrVert  ^2 }g_j(\utheta_j)\, d\utheta_j}
\\
&\leq& \int_{A_n^c}\prod_{j=0}^{\widetilde J_n(\beta)}
\nu _{j,n}^{-1}e^{G\llvert  I_j\rrvert  } \frac{ e^{-\sklfrac{ n}{2} \llVert   \utheta_j - \uY_j\rrVert
^2 }(1+\nu_{j,n})\,d\widetilde\pi_j(\utheta_j) }{ \int_{[-L_0,L_0]^{\llvert  I_j\rrvert  }} e^{-\sklfrac{ n}{2} \llVert   \utheta_j - \uY_j\rrVert  ^2
}g_j(\utheta_j)\, d\utheta_j}.
\end{eqnarray*}
On $B_n$, with obvious notation,
\begin{eqnarray*}
&&   \int_{[-L_0,L_0]^{\llvert  I_j\rrvert  }} e^{-\sklfrac{ n}{2} \llVert   \utheta_j - \uY
_j\rrVert  ^2 }d\utheta_j
\\
&&\qquad \geq  \biggl( \frac{2\pi}{n } \biggr)^{\llvert  I_j\rrvert  /2} - \operatorname{Pr} \bigl(
\exists j \leq\widetilde J_n(\beta), k \in I_j: \llvert
Y_{j,k}\rrvert + n^{-1/2} \bigl\llvert \mathcal N(0,1)\bigr
\rrvert > L_0 \bigr).
\end{eqnarray*}
Since $\llvert  \theta_{j,k}^0\rrvert   \leq L$ for all $Y^n \in B_n$, we find $\llvert
Y_{j,k}\rrvert   + n^{-1/2}\llvert  \mathcal N(0,1)\rrvert   \leq L+\frac{1}{2}+
n^{-1/2}\llvert  \mathcal N(0,1)\rrvert  $, and hence
\begin{eqnarray*}
&&   \operatorname{Pr} \bigl( \exists j \leq\widetilde J_n(\beta), k \in
I_j: \llvert Y_{j,k}\rrvert + n^{-1/2}\bigl\llvert
\mathcal N(0,1)\bigr\rrvert > L_0 \bigr)
\\
&&\qquad \leq   n\operatorname{Pr} \biggl( \bigl\llvert \mathcal N(0,1)\bigr\rrvert >
\frac{\sqrt
{n}}{2} \biggr) \leq2ne^{-n/8}.
\end{eqnarray*}
It follows that for $Y^n \in B_n$,
\[
\int_{[-L_0,L_0]^{\llvert  I_j\rrvert  }} e^{-\sklfrac{ n }{2}\llVert   \utheta_j - \uY_j\rrVert  ^2
}d\utheta_j\geq
\biggl( \frac{2\pi}{n } \biggr)^{\llvert  I_j\rrvert  /2} - 2ne^{-n/8} \geq
\frac{1}{2} \biggl( \frac{2\pi}{n } \biggr)^{\llvert  I_j\rrvert  /2}.
\]
We now study the numerator in $P^\pi( A_n^c\mid  Y^n )$. For $Y^n \in B_n$
\begin{eqnarray*}
-\llVert \utheta_j - \uY_j\rrVert ^2 &\leq&
\bigl\llVert \utheta_j^0 - \uY_j\bigr\rrVert
^2-\frac{1}2\bigl\llVert \utheta_j - \utheta
_j^0\bigr\rrVert ^2 \leq\frac{e}{n}
\llvert I_j\rrvert -\frac{1}2\bigl\llVert
\utheta_j - \utheta_j^0\bigr\rrVert
^2.
\end{eqnarray*}
On $B_n$, we can subsequently bound $P^\pi( A_n^c\mid  Y^n )$ by $2e^{\sum
_{j=0}^{\widetilde J_n(\beta)}(G+e/2)\llvert  I_j\rrvert  }$ times
\begin{eqnarray*}
&& \int_{A_n^c} \prod
_{j=0}^{\widetilde J_n(\beta)} \biggl( \frac{n}{2\pi}
\biggr)^{\llvert  I_j\rrvert  /2}\nu_{j,n}^{-1} e^{-\sklfrac{n}4 \llVert   \utheta_j - \utheta_j^0\rrVert  ^2 }(1+
\nu_{j,n})\,d\widetilde \pi_j(\utheta_j)
\\
&&\qquad = e^{-\vfrac{ nM^2u_n(\beta)^2}{32} }
\\
&&\quad\qquad{}\times \prod_{j=0}^{\widetilde J_n(\beta
)}
\int_{\R^{\llvert  I_j\rrvert  }} \biggl( \frac{n}{2\pi} \biggr)^{\llvert  I_j\rrvert  /2}
\nu _{j,n}^{-1} e^{-\sklfrac{n}8 \llVert   \utheta_j - \utheta_j^0\rrVert  ^2 }(1+\nu_{j,n})\,d
\widetilde \pi_j(\utheta_j)
\\
&&\qquad \leq e^{-\vfrac{ nM^2u_n(\beta)^2}{32} }\prod_{j=0}^{\widetilde
J_n(\beta)}
\biggl( \nu_{j,n}^{-1} \biggl( \frac{n}{2\pi}
\biggr)^{\llvert  I_j\rrvert  /2} + 2^{\llvert  I_j\rrvert  }e^{G\llvert  I_j\rrvert  } \biggr)
\\
&&\qquad \leq e^{-\vfrac{ nM^2u_n(\beta)^2}{32} }\prod_{j=0}^{\widetilde
J_n(\beta)}
\bigl( e^{c\llvert  I_j\rrvert   } + 2^{\llvert  I_j\rrvert  }e^{G\llvert  I_j\rrvert  } \bigr).
\end{eqnarray*}
Choosing $M$ large enough and using the exponential bound on $P_{\theta
_0}^n(B_n^c)$ shows that $E_{\theta_0}^n[P^{\pi}(A_n^c\mid Y^n )]\leq
e^{-An^{1/(2\beta+1)}}$. This completes the proof of Theorem~\ref{thWNL2}.

\begin{appendix}
\section*{Appendix: Additional proofs}

\subsection{Explicit bounds on \texorpdfstring{$\Omega_n$}{Omegan}} \label{metricentropy}

\mbox{}
\begin{pf*}{Proof of (\ref{OmeganforHoelder})}
Since we are on a H\"older space, we can prove the result for $\ell
=\ell_\infty$ (see also Section~\ref{secWN}). Consider $\theta
=(\theta_{j,k})_{(j,k)\in\Lambda}\in\mathcal H(\beta,L)$ and pick
$J_n(\beta)$ such that
\[
\tfrac{1}2 (L/2 M)^{1/\beta}(n/\log n)^{1/(2\beta+1)}
\leq2^{J_n(\beta
)}\leq(L/2 M)^{1/\beta}(n/\log n)^{1/(2\beta+1)}.
\]
On\vspace*{1pt} resolution level $J_n(\beta)$ chose an arbitrary index in $\Lambda
$, $(J_n(\beta),k^*)$ say. By definition of $\mathcal H(\beta,L)$
there exists $\theta'\in\mathcal H(\beta,L)$, with $\llvert  \theta
_{j,k}'-\theta_{j,k}\rrvert  $ equals $L2^{-J_n(\beta) (\beta+1/2)}$ if
$(j,k)=(J_n(\beta),k^*)$ and zero otherwise. Then $\ell_\infty
(\theta,\theta')=L2^{-J_n(\beta) \beta}\geq2\epsilon_n(\theta)$
and $\llVert  \theta-\theta'\rrVert  _{L^2}=L2^{-J_n(\beta)(\beta+1/2)}\lesssim
\sqrt{\log n/n}$.
\end{pf*}

\subsection{Derivation of \texorpdfstring{(\protect\ref{equpperboundheuristic})}{(5.6)}}
\label{subsecderivationofinequalofggvdv}
We briefly recall the main arguments of Ghosal, Ghosh and van~der Vaart
\cite{ghosalghoshvdv00}
leading to inequality (\ref{equpperboundheuristic}).
Their method is based on two assumptions, namely a bound on the local
entropy as well as existence of a decomposition $\Theta=\Theta_n\cup
(\Theta\setminus\Theta_n)$ such that the prior is uniform on $\Theta
_n$ (with respect to Kullback--Leibler balls) and assigns negligible
mass to $\Theta\setminus\Theta_n$ (cf. \cite{ghosalghoshvdv00},
equations (2.7), (2.3) and (2.5)). To derive (\ref{equpperboundheuristic}) only the assumption on the prior needs to be imposed.
Recall that
\[
P^\pi \bigl(\theta: \ell( \theta_0, \theta) >
\epsilon_n \mid Y^n \bigr) =\frac{ \int_{\ell( \theta_0, \theta) > \epsilon_n} e^{ \mathcal
L_n(\theta) - \mathcal L_n(\theta_0) }\pi(d\theta)}{ \int_{\Theta
} e^{ \mathcal L_n(\theta) - \mathcal L_n(\theta_0) }\pi(d\theta)} =:
\frac{N_n }{ D_n}.
\]
Under the imposed conditions, there are constants $c,c'>0$ such that
$P_{\theta_0}^n(D_n\geq\exp(-cnu_n^2))\geq1-e^{-c'nu_n^2}$ (cf.
\cite{ghosalghoshvdv00}, Lemma 8.4). Hence, for any test function~$\Phi_n$,
\begin{eqnarray*}
&& E_{\theta_0}^n \bigl[P^\pi \bigl(\theta: \ell(
\theta_0, \theta) > \epsilon_n \mid Y^n \bigr)
\bigr]
\\
&&\qquad  \leq E_{\theta_0}^n[\Phi_n] + e^{cnu_n^2}
E_{\theta_0}^n \biggl[\int_{\ell(\theta_0,\theta) >\epsilon_n}
e^{ \mathcal L_n(\theta
) - \mathcal L_n(\theta_0) }(1 - \Phi_n) \pi(d\theta) \biggr]+e^{-c'nu_n^2}
\\
&&\qquad  \leq E_{\theta_0}^n[\Phi_n]+e^{cnu_n^2} \sup
_{\theta: \ell(
\theta_0, \theta) > \epsilon_n } E_\theta^n [ 1 - \Phi_n
]+e^{-c'nu_n^2}.
\end{eqnarray*}


\subsection{Proof of Lemma \texorpdfstring{\protect\ref{lemWNSc}}{1}} \label{appWNSL}

%

\mbox{}

\begin{pf*}{Proof of (\ref{WNSc1})}
%
We have
\begin{eqnarray*}
P^\pi \bigl( S^c \cap\mathcal J_n(\overline
\gamma) \neq\varnothing \mid Y^n \bigr) &\leq&\sum
_{(j,k)\in\mathcal J_n(\overline\gamma
)} P^\pi\bigl( \theta_{j,k}=0 \mid
Y^n\bigr)
\\
&\leq& \sum_{(j,k)\in\mathcal J_n(\overline\gamma)} \frac{ e^{ -\sklfrac{
n}{2} Y_{j,k}^2 }}{ w_{j,n} \int_{\R} e^{-\sklfrac{ n}{2} ( Y_{j,k} -
\theta)^2 }g(\theta)\,d\theta}.
\end{eqnarray*}
Recall\vspace*{1pt} the definition of $\Omega_{n,B}$ in the proof of Theorem~\ref
{thWNlinf}. If $Y^n\in\Omega_{n,B}$ and $\overline\gamma$ is
large enough, then $\llvert  Y_{j,k}\rrvert  >\frac{1}2 \overline\gamma\sqrt{\log
n/n}$. With the same argument as in (\ref{lowerbounddenominator}),
\[
P^\pi \bigl( S^c \cap\mathcal J_n(\overline
\gamma) \neq\varnothing \mid Y^n \bigr) \leq \sum
_{j\leq J_n, k \in I_j} \frac{ e^{ -\vfrac{ n Y_{j,k}^2 }{2}}
\sqrt{n}}{ w_{j,n}a \sqrt{\pi} } \leq\frac{ n^{K+3/2-\overline
\gamma^2/8} }{a\sqrt{\pi}},
\]
and together with (\ref{Omegan}) this completes the proof of (\ref
{WNSc1}), provided $\overline\gamma$ is sufficiently large.
\end{pf*}

%

\begin{pf*}{Proof of (\ref{WNSc2})}
We have
\begin{eqnarray*}
P^\pi \bigl( S \cap\mathcal J_n(\underline
\gamma)^c \neq \varnothing\mid Y^n \bigr) &=& \sum
_{(j,k)\in\mathcal J_n(\underline\gamma)^c} P^\pi\bigl( \theta _{j,k} \neq0
\mid Y^n\bigr)
\\
&\leq&\sum_{(j,k) \in\mathcal J_n(\underline\gamma)^c} \frac{
w_{j,n} \int_{\R} e^{ - \sklfrac{n}{2} ( \theta- Y_{j,k})^2} g(\theta
) \,d\theta}{ (1 -w_{j,n}) e^{ - \sklfrac{n}{2} Y_{j,k}^2} }
\\
&\leq& 2\sqrt{2\pi}n^{-1/2}\sup_xg(x) \sum
_{(j,k) \in\mathcal
J_n(\underline\gamma)^c } w_{j,n} e^{\vfrac{nY_{j,k}^2}{2} }.
\end{eqnarray*}
If $Y^n \in\Omega_{n,B}$, for any $(j,k) \in\mathcal J_n(\underline
\gamma)^c$,
\begin{eqnarray*}
nY_{j,k}^2 &\leq&\underline\gamma^2 \log n + 2
\log\llvert I_j\rrvert + 2B\log n + 2 \underline\gamma\sqrt{\log n}
\sqrt{2\log\llvert I_j\rrvert + 2B\log n}
\\
& \leq&2\log\llvert I_j\rrvert + \bigl( 2B + \underline
\gamma^2+ \underline\gamma C\bigr)\log n,
\end{eqnarray*}
for some constant $C$. Hence, whenever $Y^n \in\Omega_{n,B}$, using
that $w_{j,n}\leq\break  n^{(\tau\wedge1)/2-1/4-2B}2^{-j(1+\tau)}$ with
$\tau>1/2$,
\begin{eqnarray*}
P^\pi \bigl( S \cap\mathcal J_n(\underline
\gamma)^c \neq \varnothing\mid Y^n \bigr) &\lesssim&
n^{-3/4+(\tau\wedge1)/2-B+\underline \gamma C/2+
\underline\gamma^2/2} \sum_{j\leq J_n} \llvert I_j
\rrvert ^2 2^{-j(1+\tau)}
\\
&\lesssim & n^{1/4-(\tau\wedge1)/2-B + \underline \gamma C/2+
\underline\gamma^2/2} = O\bigl(n^{-B}\bigr),
\end{eqnarray*}
where for the last equality, we need that $\underline\gamma$ is
sufficiently small. The proof of (\ref{WNSc2}) follows from (\ref{Omegan}).
\end{pf*}

\subsection{Proof of Proposition \texorpdfstring{\protect\ref{upbmodel}}{1}} \label{appprlemupb}
We start with verifying condition (\ref{sievecond2}). For \mbox{$r\geq1$},
there exists an injective mapping
$\psi: \mathcal I_r \rightarrow\mathcal I_0$ such that
\[
\psi(\theta)_{\mathcal U} = \theta^*_{\mathcal U}\quad\mbox {and}\quad
\bigl\llvert \psi(\theta)_{j,k} -\theta_{j,k}\bigr\rrvert \neq
\phi_n\qquad \forall(j,k) \notin\mathcal U.
\]
%
This implies in particular that $\llvert  \mathcal J_r \rrvert   \leq\llvert  \mathcal
I_r\rrvert  \leq\llvert  \mathcal I_0\rrvert   \leq\llvert  \mathcal J_0\rrvert  $ and the partition is
admissible. Also, by construction of $\mathcal D_n$, for every $\theta
\in\mathcal J_r$ we have that $\ell_2(\theta, \psi(\theta))^2$
takes its values in the lattice $\{\phi_n^2, 2\phi_n^2, 3\phi
_n^2,\ldots\}$ and the cardinality of $\{r:u_r^2 = \phi_n^2\}$ is
bounded by $2\sum_{j\leq J_n}\llvert  I_j\rrvert   = I$. By induction on $M
=1,2,\ldots,$ the cardinality of $\{ r: u_r^2 = M\phi_n^2\}$ is
further bounded by $\sum_{i=1}^M I^i\leq(Cn)^{M+1} $ for some $C>0$.
This implies that for any $K_0>0$,
\[
\sum_{r=1}^{R_n} e^{- K_0n u_r^2} \leq\sum
_{M\geq1} e^{- K_0nM\phi
_n^2} \bigl\llvert \bigl\{ r:
u_r^2 = M \phi_n^2\bigr\}\bigr
\rrvert \leq\sum_{M\geq
1}n^{-K_0\phi_0M}
(Cn)^{M+1},
\]
which has polynomial decay in $n$ as soon as $\phi_0>2/K_0$, and can
thus be taken of the form $e^{- K_1 n \Omega(\epsilon_n(\cdot),
\mathcal H(\beta, L), \ell_\infty)^2 }$ for some $K_1>0$.
This bound is not based on any specific assumption on the experiment $\{
P_\theta^n,\theta\in\Theta\}$ and only depends on the set $\Theta
$, the loss $\ell= \ell_\infty$, and $d = \ell_2$. It remains to
check condition (\ref{sievecond1}). We first consider the white noise
model. Then
\begin{eqnarray*}
&&  -n^{-1} \bigl(\mathcal L_n(\theta) - \mathcal
L_n \bigl(\psi(\theta ) \bigr) \bigr)
\\
&&\qquad =   \frac{\llVert   \theta - \theta_0\rrVert  _{L^2}^2 - \llVert   \psi(\theta) -
\theta_0\rrVert  _{L^2}^2}{2} - \sum_{(j,k) \in\Lambda} \bigl(
Y_{j,k}-\theta ^0_{j,k}\bigr) \bigl(
\theta_{j,k} - \psi(\theta)_{j,k}\bigr)
\\
&&\qquad = \frac{\llVert   \theta -\psi( \theta)\rrVert  _{L^2}^2}{2} +\bigl\langle\theta- \psi(\theta), \psi(\theta) -
\theta_0\bigr\rangle_{L^2}
\\
&&\quad\qquad{} - \sum
_{(j,k)
\in\Lambda} \bigl( Y_{j,k}-\theta^0_{j,k}
\bigr) \bigl( \theta_{j,k} - \psi(\theta)_{j,k}\bigr).
\end{eqnarray*}
The above computation is simply a sequential formulation of the
Cameron--Martin formula: Here, we emphasise on the property that $\ell
_2(\theta, \theta') = \llVert  \theta-\theta'\rrVert  _{L^2}$ is a Hilbert norm
associated to the scalar product $\langle\cdot,\cdot\rangle_{L^2}$.
The sum in $(j,k)\in\Lambda$ involving $Y_{j,k}$ has to be understood
as a limit in $L^2(P_\theta^n)$, and it is well defined since $\theta
-\psi(\theta)\in\ell^2(\Lambda)$ and the $Y_{j,k}$ are independent
and standard normal under $P_\theta^n$.

Recall\vspace*{1pt} the definition of $\mathcal U$ in (\ref{defUinconstrofsieve}). For $(j,k)\in\mathcal U$, we have by construction $\llvert  \theta
_{j,k}^0-\psi(\theta)_{j,k}\rrvert  \leq\phi_n/4$ and for $(j,k)\in
\mathcal U^c$, $\llvert  \theta_{j,k}^0-\psi(\theta)_{j,k}\rrvert  \leq3\phi_n/4$.
In the latter case, we also know that $\llvert  \theta_{j,k}-\psi(\theta
)_{j,k}\rrvert  \neq\phi_n$ but has values in $\{0,2\phi_n,3\phi_n,\ldots\}
$. Therefore,
\begin{eqnarray*}
\mathcal L_n(\theta) - \mathcal L_n\bigl(\psi(\theta)
\bigr) & \leq& - \frac{ n\llVert   \theta -\psi( \theta)\rrVert  _{L^2}^2}{ 8 } + n \sum_{(j,k) \in\Lambda}
\bigl( Y_{j,k}-\theta^0_{j,k}\bigr) \bigl(
\theta_{j,k} - \psi (\theta)_{j,k}\bigr).
\end{eqnarray*}
Introduce the event $\Omega_n = \{\max_{j\leq J_n, k\in I_j }
\llvert  Y_{j,k} -\theta^0_{j,k}\rrvert   \sqrt{n}\leq2 \sqrt{\log n} \}$. For
$Y^n\in\Omega_n$,
\begin{eqnarray*}
&&  \biggl\llvert \sum_{(j,k) \in\Lambda} \bigl(
Y_{j,k}-\theta^0_{j,k}\bigr) \bigl( \theta
_{j,k} - \psi(\theta)_{j,k}\bigr) \biggr\rrvert
\\
&&\qquad = \biggl\llvert \sum_{(j,k) \in\Lambda} {\mathbf1}_{ \llvert   \theta_{j,k} - \psi
(\theta)_{j,k}\rrvert   \geq\phi_n}
\bigl( Y_{j,k}-\theta^0_{j,k}\bigr) \bigl( \theta
_{j,k} - \psi(\theta)_{j,k}\bigr)\biggr\rrvert
\\
&&\qquad \leq  2\phi_0^{-1}\bigl\llVert \theta- \psi(\theta)\bigr
\rrVert _{L^2}^2
\end{eqnarray*}
due to $\llvert  Y_{j,k}-\theta^0_{j,k}\rrvert   \leq2\sqrt{\log n}/\sqrt{n} \leq
2/\phi_0 \phi_n$. Picking $\phi_0$ large enough,
%
\[
\mathcal L_n(\theta) - \mathcal L_n \bigl(\psi(\theta)
\bigr) \leq- \frac{ n\llVert   \theta -\psi( \theta)\rrVert  _{L^2}^2}{ 16 }\qquad\mbox{on }\Omega_n.
\]
Since $P_{\theta_0}^n(\Omega_n^c)\leq2n^{-1}$, condition (\ref
{sievecond1}) is satisfied. This completes the proof of Proposition
\ref{upbmodel}.
%

\subsection{Proof of Proposition \texorpdfstring{\protect\ref{teodensityestimation}}{3}} \label{appthmthmupbdensityestimation}

We start with the construction of the prior $\pi_n$. Contrariwise to
the white noise model, we truncate $j\leq\overline J_n$ with $\sqrt
{n}/\log n< 2^{\overline J_n}\leq2\sqrt{n}/\log n$. Set $\Theta=
\bigcup_{\beta\in[\beta_1, \beta_2]} \mathcal H'(\beta, L)$.
Recall the definition of $\mathcal D_n$ in (\ref{Dnsievedef}) with
$\phi_n=\phi_0\sqrt{\log n/n}$ and consider
\[
\mathcal D_n'=\bigl\{\theta\in\mathcal
D_n: \exists\theta'\in\Theta \mbox{ such that }  \bigl\llvert \theta_{j,k}-\theta_{j,k}'\bigr
\rrvert \leq\phi_n, \forall j\leq\overline J_n, k\in
I_j\bigr\}
\]
as set of nonnormalised test densities. By construction $\sum_{j,k}
\theta_{j,k}\Psi_{j,k}\geq c/2$, $\forall\theta\in\mathcal D_n'$
and, therefore, $\sqrt{f_\theta}= \llVert  \theta\rrVert  _{L^2}^{-1} \sum_{j,k}
\theta_{j,k}\Psi_{j,k}$ is well-defined, that is, $f_\theta$ is a
density [note that this definition extends (\ref{eqsqrtdensexpansion}) in a consistent way]. The set $\mathcal D_n'$ constitutes
the sieve and the prior is given by $\pi_n \propto\sum_{\theta\in
\mathcal D_n'} \delta_{\theta/\llVert  \theta\rrVert  _{L^2}}$.

For the subsequent analysis, we need some inequalities for the elements
in $\mathcal D_n'$, which are derived next. Due to $\beta_1>1/2$, the
coefficients of the parameter vectors are absolutely summable and
\[
\overline L = \max \biggl( \sup_{\theta\in\Theta\cup\mathcal D_n'} \ell_\infty
\bigl(\theta/\llVert \theta\rrVert _{L^2},0\bigr)+\sum
_{(j,k)\in\Lambda} \llvert \theta_{j,k}\rrvert, 1 \biggr) <
\infty.
\]
Let $\theta\in\mathcal D_n'$. By construction, there exists a $\theta
'\in\Theta$ with $\llVert  \theta'\rrVert  _{L^2}=1$ and $\llvert  \theta_{j,k}-\theta
_{j,k}'\rrvert  \leq\phi_n$ for all $(j,k)\in\Lambda$. With $\llVert  \theta\rrVert
_{L^2}^2=\langle\theta+\theta',\theta-\theta'\rangle_{L^2}+1$ we find
%
\begin{equation}
\label{eqthetaL^2ineqs}
\qquad\tfrac{1}2 \leq\llVert \theta\rrVert _{L^2}
\leq2,\qquad\bigl\llvert \llVert \theta\rrVert _{L^2}-1\bigr\rrvert
\leq4\overline L \phi_n\quad\mbox{and}\quad\biggl\llvert
\frac{1}{\llVert  \theta\rrVert  _{L^2}}-1\biggr\rrvert \leq8\overline L\phi_n.
\end{equation}
Next, let us construct an admissible partition. Notice that there is a
finite $J_0$, such that
%
\begin{equation}
\label{eqdefofJ0} \sup_{\theta\in\Theta\cup\mathcal D_n'} \max_{j>J_0, k\in I_j}
\llvert \theta_{j,k}\rrvert + \sum_{j>J_0}^\infty
\sum_{k\in I_j} \theta _{j,k}^2<2^{-21}
\frac{1}{\overline L^3}.
\end{equation}
Let $Q=\lceil\overline L^22^{11}\rceil$. For every $(j,k)$, we can
define an equivalence relation $\simeq$ via $\theta_{j,k}\simeq
\theta_{j,k}'$ iff $\theta_{j,k}=\theta_{j,k}'$ or $\theta
_{j,k},\theta_{j,k}'\in(\theta_{j,k}^0-q_{j,k}(\theta_0)\phi_n,
\theta_{j,k}^0+q_{j,k}(\theta_0)\phi_n)$, where
\[
q_{j,k}(\theta_0)= \cases{ Q, &\quad if $j\leq
J_0$,
\cr
1, &\quad if $j> J_0$, $\bigl\llvert
\theta_{j,k}^0\bigr\rrvert >2^{-9}
\phi_n$,
\cr
0, &\quad if $j> J_0$, $\bigl\llvert
\theta_{j,k}^0\bigr\rrvert \leq2^{-9}
\phi_n$.}
\]
This induces an equivalence relation on the nonnormalised densities
$\mathcal D_n'$ via $\theta\simeq\theta'$ iff $\theta_{j,k}\simeq
\theta_{j,k}'$ for all $(j,k)$, $j\leq\overline J_n$. By
construction, there exists $\theta^*\in\mathcal D_n'$ such that
$\llvert  \theta_{j,k}^0-\theta_{j,k}^*\rrvert  \leq\frac{1}2 \phi_n$ for all
$(j,k)$. Denote by $\mathcal I_r$, $r= 0, 1, \ldots$ the equivalence
classes of $\mathcal D_n'$ and let $\mathcal I_0$ be the equivalence
class of $\theta^*$. Define $J_n(\beta)$ as in (\ref{eqJndefproperties}), replacing $\phi_n/4$ by $2^{-9}\phi_n$ in the first
condition. Using (\ref{eqthetaL^2ineqs}), there exists a constant
$A=A(\beta, \overline L, \phi_0, Q)$ such that, for all $\theta\in
\mathcal I_0$,
\begin{eqnarray*}
\ell_{\infty}\bigl(\theta_0,\theta/\llVert \theta\rrVert
_{L^2}\bigr) &\leq&\ell_{\infty}(\theta_0,\theta)+4
\overline L \phi_n \ell _\infty\bigl(\theta/\llVert \theta
\rrVert _{L^2},0\bigr)
\\
&\leq& 2Q\phi_n \sum_{j=0}^{J_n(\beta)}
2^{j/2} +\epsilon_n(\beta)+ 4\overline L^2
\phi_n
\\
&\leq& A\epsilon_n(\beta).
\end{eqnarray*}
For this $A$, we define $\mathcal J_0 = \{ \theta\in\mathcal D_n':
\ell_\infty(\theta_0, \theta/\llVert  \theta\rrVert  _{L^2}) \leq A\epsilon
_n(\beta)\}$, $\mathcal J_r = \mathcal I_r \cap\mathcal J_0^c$. Now,
for any $r\geq1$, we construct an injective map $\psi: \mathcal
J_r\rightarrow\mathcal J_0$ and verify that for this map the
properties (\ref{sievecond1}) and (\ref{sievecond2}) hold. To this
end, define $\iota(\theta_{j,k})$ as $\lceil\theta_{j,k}^0 \phi
_n^{-1}\rceil\phi_n$ if $\theta_{j,k}>\theta_{j,k}^0$ and $\lfloor
\theta_{j,k}^0 \phi_n^{-1}\rfloor\phi_n$ otherwise. If $(j,k) \in
\mathcal J_r$, $r\neq0$,
\[
\psi(\theta)_{j,k}= \cases{ \theta_{j,k}, &\quad if $\bigl
\llvert \theta_{j,k}-\theta _{j,k}^0\bigr\rrvert
<q_{j,k}(\theta_0)\phi_n$,
\vspace*{3pt}\cr
\iota(
\theta_{j,k}), &\quad if $ \bigl\llvert \theta_{j,k}-
\theta_{j,k}^0\bigr\rrvert \geq q_{j,k}(
\theta_0)\phi_n,  q_{j,k}(\theta_0)>0$,
\vspace*{3pt}\cr
0, &\quad if $q_{j,k}(\theta_0)=0$.} %
\]
It is not difficult to see that $\psi: \mathcal J_r\rightarrow
\mathcal J_0$ is injective. This completes the proof of the admissible
part. By following the same arguments as in the proof of Proposition~\ref{upbmodel}, condition (\ref{sievecond2}) can be verified.

Therefore, it remains to check (\ref{sievecond1}). For $u>0$, we have
$\log u=2\log(\sqrt{u})\leq2(\sqrt{u}-1)$ and, therefore,
%
\begin{equation}
\label{eqfirstupperbdLRdensest} \mathcal L_n(\theta)-\mathcal L_n\bigl(
\psi(\theta)\bigr)\leq2 \sum_{i=1}^n
\biggl(\frac{\sqrt{f_\theta}(Y_i)}{\sqrt{f_{\psi(\theta
)}}(Y_i)}-1 \biggr).
\end{equation}
We further decompose the right-hand side using
%
\begin{eqnarray}\label{eqLRdecomp}
\frac{x}y-1 &=& \frac{x-y}z+\frac{(x-y)(z-y)}{z^2}
\nonumber\\[-8pt]\\[-8pt]\nonumber
&&{} + \frac {(x-y)(z-y)^2}{z^3}+\frac{(x-y)(z-y)^3}{z^3y}
\end{eqnarray}
with $x=\sqrt{f_\theta}(Y_i)$, $y=\sqrt{f_{\psi(\theta)}}(Y_i)$,
and $z=\sqrt{f_{\theta_0}}(Y_i)$. In the sequel, we control the large
deviations behaviour of the terms on the right-hand side separately
[denoting the single steps by (I)--(IV)]. The key ingredient is the
following well-known version of Bernstein's inequality: If $X_1,\ldots,X_n$ are i.i.d., centered and $\llvert  X_i\rrvert  \leq M$, then $\forall t>0$,
$\mathbb P(\llvert  \sum_{i=1}^n X_i\rrvert  >t)\leq2\exp(-\frac{1}2 t^2/\break (n\mathbb
E[X_1^2]+Mt/3))$.

\begin{longlist}[(III)]
\item[(I)] Define $\Omega_{n,1}(\tau)$ as the event
\[
\Biggl\{\Biggl\llvert \sum_{i=1}^n
\frac{\Psi_{j,k}(Y_i)}{\sqrt{f_{\theta
_0}}(Y_i)}-n \int\Psi_{j,k}(u) \sqrt{f_{\theta_0}}(u)\,du
\Biggr\rrvert \leq \tau\sqrt{n\log n}, \forall j\leq\overline J_n, k
\in I_j\Biggr\}.
\]
Observe\vspace*{1pt} that the random variables $\Psi_{j,k}(Y_i)/\sqrt{f_{\theta
_0}}(Y_i)$, $i=1,\ldots,n$, are i.i.d., bounded in absolute value by a
multiple of $n^{1/4}$ and their second moment is one. Thus, by a union
bound and Bernstein's inequality, $P_{\theta_0}^n(\Omega_{n,1}(\tau
)^c)\lesssim n^{-1}$, provided that $\tau$ is large enough. On $Y^n\in
\Omega_{n,1}(\tau)$,
\begin{eqnarray*}
\sum_{i=1}^n \frac{\sqrt{f_\theta} (Y_i)-\sqrt{f_{\psi(\theta
)}}(Y_i)}{\sqrt{f_{\theta_0}}(Y_i)} &\leq&n
\int \bigl(\sqrt{f_\theta}(u)-\sqrt{f_{\psi(\theta
)}}(u) \bigr)
\sqrt{f_{\theta_0}}(u) \,du
\\
&&{}+\tau\sqrt{n \log n} \sum_{(j,k)} \biggl\llvert
\frac{\theta_{j,k}}{\llVert
\theta\rrVert  _{L^2}}-\frac{\psi(\theta)_{j,k}}{\llVert  \psi(\theta)\rrVert
_{L^2}}\biggr\rrvert.
\end{eqnarray*}
Using the inequalities (\ref{eqthetaL^2ineqs}), we can bound the
second term on the right-hand side by $\tau\phi_0^{-1}(1+16\overline
L^2)n \llVert  \theta-\psi(\theta)\rrVert  _{L^2}^2$, and thus, making $\phi_0$
large enough we obtain on $Y^n \in\Omega_{n,1}(\tau)$,
\begin{eqnarray*}
&&\sum_{i=1}^n \frac{\sqrt{f_\theta} (Y_i)-\sqrt{f_{\psi(\theta
)}}(Y_i)}{\sqrt{f_{\theta_0}}(Y_i)}
\\
&&\qquad \leq n \int \bigl(\sqrt{f_\theta}(u)-\sqrt{f_{\psi(\theta
)}}(u) \bigr)
\sqrt{f_{\theta_0}}(u) \,du+2^{-9}n \bigl\llVert \theta-\psi(\theta
)\bigr\rrVert _{L^2}^2.
\end{eqnarray*}

\item[(II)] Similar as for (I), set $\Omega_{n,2}$ for the event
\begin{eqnarray*}
&& \Biggl\{\Biggl\llvert \sum_{i=1}^n
\frac{\Psi_{j,k}(Y_i)\Psi
_{j',k'}(Y_i)}{f_{\theta_0}(Y_i)}-n \delta_{(j,k),(j',k')}\Biggr\rrvert \leq n^{3/4}
\log n,
\\
&&\qquad \forall j,j'\leq\overline J_n, k\in
I_j, k'\in I_{j'}\Biggr\},
\end{eqnarray*}
with $\delta_{(j,k),(j',k')}$ the Kronecker delta. Now, $\Psi
_{j,k}(Y_i)\Psi_{j',k'}(Y_i)/f_{\theta_0}(Y_i)$, $i=1,\ldots,n$, are
i.i.d. and bounded in absolute value by a multiple of $\sqrt{n}$. The
second moment is also smaller than $\mathrm{const.}\times n^{1/2}$. Using a
union bound and Bernstein's inequality, $P_{\theta_0}^n(\Omega
_{n,2}^c)\lesssim n^{-1}$ for $n$ large enough. On $Y^n \in\Omega
_{n,2}$, we see that $\sum_{i=1}^n (\sqrt{f_\theta}(Y_i)-\sqrt
{f_{\psi(\theta)}}(Y_i))(\sqrt{f_{\theta_0}}(Y_i)-\sqrt{f_{\psi
(\theta)}}(Y_i))/f_{\theta_0}(Y_i)$ can be bounded by its expectation plus
\[
\sum_{j,k} \biggl\llvert \frac{\theta_{j,k}}{\llVert  \theta\rrVert  _{L^2}}-
\frac{\psi
(\theta)_{j,k}}{\llVert  \psi(\theta)\rrVert  _{L^2}}\biggr\rrvert \sum_{j',k': \Psi_{j,k}\Psi_{j',k'}\neq0} \biggl
\llvert \theta_{j,k}^0- \frac{\psi(\theta)_{j,k}}{\llVert  \psi(\theta)\rrVert  _{L^2}}\biggr\rrvert
n^{3/4}\log n.
\]
Due to the compact support of $\Psi$, there are of the order of $\log
n$ index pairs $(j',k')$ with $j'\leq\overline J_n$ and $\Psi
_{j,k}\Psi_{j',k'}\neq0$. Using that $\psi(\theta)\in\mathcal
J_0$, together with the inequalities (\ref{eqthetaL^2ineqs}),
yields $\llvert  \theta_{j,k}^0-\psi(\theta)_{j,k}/\llVert  \psi(\theta)\rrVert
_{L^2}\rrvert  \leq(Q+4\overline L)\phi_n$. Because of $\llvert  \theta_{j,k}-\psi
(\theta)_{j,k}\rrvert   {\mathbf1}_{\theta_{j,k}\neq\psi(\theta)_{j,k}}\geq
\phi_n$, the expression in the last display is smaller than $2^{-9}n \llVert
\theta-\psi(\theta)\rrVert  _{L^2}^2$, for sufficiently large $n$ and,
therefore, on $Y^n \in\Omega_{n,2}$,
\begin{eqnarray*}
&& \sum_{i=1}^n \frac{(\sqrt{f_\theta}(Y_i)-\sqrt{f_{\psi(\theta
)}}(Y_i))(\sqrt{f_{\theta_0}}(Y_i)-\sqrt{f_{\psi(\theta
)}}(Y_i))}{f_{\theta_0}(Y_i)}
\\
&&\qquad \leq n\int\bigl(\sqrt{f_\theta}(u)-\sqrt{f_{\psi(\theta)}}(u)\bigr)
\\
&&\quad\qquad{}\times
\bigl(\sqrt {f_{\theta_0}}(u)-\sqrt{f_{\psi(\theta)}}(u)\bigr)\,du
+2^{-9}n\bigl\llVert \theta-\psi(\theta)\bigr\rrVert
_{L^2}^2.
\end{eqnarray*}

\item[(III)] This case works similar as (II) and is therefore only
sketched here. In fact, we need to consider $\Omega_{n,3}$ which is
the same event as $\Omega_{n,2}$ but applied to the random variables
$\Psi_{j_1,k_1}(Y_i)\Psi_{j_2,k_2}(Y_i)\Psi_{j_3,k_3}(Y_i)/f_{\theta
_0}^{3/2}(Y_i)$ (and the $n^{3/4}$ should be exchanged with $n$). Since
these random variables are bounded in absolute value by a constant
times $n^{3/4}$ and have second moment smaller than a constant times
$n$, we obtain $P_{\theta_0}^n (\Omega_{n,3}^c)\lesssim n^{-1}$.
Using the inequalities (\ref{eqthetaL^2ineqs}) again, on $Y^n \in
\Omega_{n,3}$,
\begin{eqnarray*}
&& \sum_{i=1}^n \frac{(\sqrt{f_\theta}(Y_i)-\sqrt{f_{\psi(\theta
)}}(Y_i))(\sqrt{f_{\theta_0}}(Y_i)-\sqrt{f_{\psi(\theta
)}}(Y_i))^2}{f_{\theta_0}^{3/2}(Y_i)}
\\
&&\qquad \leq n\int\bigl(\sqrt{f_\theta}(u)-\sqrt{f_{\psi(\theta)}}(u)\bigr)
\bigl(\sqrt {f_{\theta_0}}(u)-\sqrt{f_{\psi(\theta)}}(u)\bigr)^2
\frac{du}{\sqrt
{f_{\theta_0}}(u)}
\\
&&\quad\qquad{}  + 2^{-10}n\bigl\llVert \theta-\psi(\theta)\bigr\rrVert
_{L^2}^2.
\end{eqnarray*}
The first term on the right-hand side can be further bounded by a
constant times
\[
n \sum_{j,k} \biggl\llvert \frac{\theta_{j,k}}{\llVert  \theta\rrVert  _{L^2}}-
\frac
{\psi(\theta)_{j,k}}{\llVert  \psi(\theta)\rrVert  _{L^2}}\biggr\rrvert \int\bigl\llvert \Psi_{j,k}(u)\bigr
\rrvert \bigl(\sqrt{f_{\theta_0}}(u)-\sqrt{f_{\psi(\theta
)}}(u)
\bigr)^2 \,du.
\]
Expanding $\sqrt{f_{\theta_0}}(u)-\sqrt{f_{\psi(\theta)}}(u)$ and
using the compactness of $\Psi$ as well as $\llvert  \theta_{j,k}^0-\psi
(\theta)_{j,k}/\llVert  \psi(\theta)\rrVert  _{L^2}\rrvert  \leq(Q+4\overline L)\phi_n$
and (\ref{eqthetaL^2ineqs}), we find that the last display can be
further bounded by $O(n\llVert  \theta-\psi(\theta)\rrVert  _{L^2}^2 \phi_n \log
n)$, and so, on $Y^n \in\Omega_{n,3}$,
\begin{eqnarray*}
&&\sum_{i=1}^n \frac{(\sqrt{f_\theta}(Y_i)-\sqrt{f_{\psi(\theta
)}}(Y_i))(\sqrt{f_{\theta_0}}(Y_i)-\sqrt{f_{\psi(\theta
)}}(Y_i))^2}{f_{\theta_0}^{3/2}(Y_i)}
\leq2^{-9}n\bigl\llVert \theta-\psi(\theta)\bigr\rrVert
_{L^2}^2.
\end{eqnarray*}

\item[(IV)] For this term, no exponential inequality is needed, and a
deterministic bound can be obtained as follows. Observe that there is a
constant $c(\Psi)$, such that $\llvert  \sqrt{f_{\theta_0}}(Y_i)-\sqrt
{f_{\psi(\theta)}}(Y_i)\rrvert  \leq c(\Psi)Q\phi_n \sum_{j=0}^{\overline
J_n} 2^{j/2}\lesssim2^{\overline J_n/2}\phi_n$. This shows that
\[
\sum_{i=1}^n \bigl(\sqrt{f_\theta}(Y_i)-
\sqrt{f_{\psi(\theta
)}}(Y_i)\bigr) \bigl(\sqrt{f_{\theta_0}}(Y_i)-
\sqrt{f_{\psi(\theta
)}}(Y_i)\bigr)^3/
\bigl(f_{\theta_0}^{3/2}(Y_i)\sqrt{f_{\psi(\theta)}}(Y_i)
\bigr)
\]
can be bounded by a constant times
\begin{eqnarray*}
&& n2^{2\overline J_n} \phi_n^3\sum
_{j,k} \biggl\llvert \frac{\theta_{j,k}}{\llVert
\theta\rrVert  _{L^2}}-\frac{\psi(\theta)_{j,k}}{\llVert  \psi(\theta)\rrVert
_{L^2}}
\biggr\rrvert.
\end{eqnarray*}
Using the definition of $\overline J_n$ and (\ref{eqthetaL^2ineqs}), we find that this term is of negligible order $O(n\llVert  \theta
-\psi(\theta)\rrVert  _2^2/\log n)$, uniformly over $\theta$.

Now, we are ready to complete the proof. Since $P_{\theta_0}^n((\Omega
_{n,1}(\tau)\cap\Omega_{n,2}\cap\Omega_{n,3})^c)\lesssim n^{-1}$,
we can throughout the following assume that $Y^n \in\Omega_{n,1}(\tau
)\cap\Omega_{n,2}\cap\Omega_{n,3}$. It is then enough to prove
$\frac{1}n  (\mathcal L_n(\theta)-\mathcal L_n(\psi(\theta))
)\leq-K_0\llVert  \theta-\psi(\theta)\rrVert  _{L^2}^2$ for\vspace*{1pt} some positive
constant $K_0$. Combining the estimates in (I)--(IV), with (\ref
{eqfirstupperbdLRdensest}) and (\ref{eqLRdecomp}), we find,
for sufficiently large $n$,
\begin{eqnarray*}
\hspace*{-3pt} && \frac{1}n \bigl(\mathcal L_n(\theta)-\mathcal
L_n\bigl(\psi(\theta)\bigr) \bigr)
\\
\hspace*{-3pt} &&\qquad \leq \int \bigl(\sqrt{f_\theta}(u)-\sqrt{f_{\psi(\theta)}}(u)\bigr)
\bigl(2\sqrt{f_{\theta
_0}}(u)-\sqrt{f_{\psi(\theta)}}(u)\bigr) \,du
+2^{-7} \bigl\llVert \theta-\psi(\theta)\bigr\rrVert
_{L^2}^2
\\
\hspace*{-3pt} &&\qquad = \biggl(\frac{1}2 -\bigl\llVert \psi(\theta)\bigr\rrVert
_{L^2}\biggr) \sum_{j,k} \biggl(
\frac{\theta_{j,k}}{\llVert  \theta\rrVert  _{L^2}}-\frac{\psi
(\theta)_{j,k}}{\llVert  \psi(\theta)\rrVert  _{L^2}} \biggr)^2
\\
\hspace*{-3pt} &&\quad\qquad {} + 2 \sum_{j,k} \biggl(\frac{\theta_{j,k}}{\llVert  \theta\rrVert  _{L^2}}-
\frac{\psi
(\theta)_{j,k}}{\llVert  \psi(\theta)\rrVert  _{L^2}} \biggr) \bigl(\theta_{j,k}^0-\psi (
\theta)_{j,k}\bigr)+2^{-7} \bigl\llVert \theta-\psi(\theta)
\bigr\rrVert _{L^2}^2
\end{eqnarray*}
using that
\[
\int\bigl(\sqrt{f_\theta}(u)-\sqrt{f_{\psi(\theta)}}(u)\bigr)
\sqrt{f_{\psi
(\theta)}}(u) \,du=-\frac{1}2 \int\bigl(\sqrt{f_\theta}(u)-
\sqrt{f_{\psi
(\theta)}}(u)\bigr)^2 \,du.
\]
If $\llvert  \theta_{j,k}^0\rrvert  >2^{-9}\phi_n$, then by construction of $\psi
(\theta)_{j,k}$, we have
\[
\bigl(\theta_{j,k}-\psi(\theta)_{j,k}\bigr) \bigl(
\theta_{j,k}^0-\psi(\theta )_{j,k}\bigr)\leq0.
\]
Otherwise, if $\llvert  \theta_{j,k}^0\rrvert  \leq2^{-9}\phi_n$, then $\psi(\theta
)_{j,k}=0$ and so
\[
\frac{2}{\llVert  \psi(\theta)\rrVert  _{L^2}}\sum_{j,k}\bigl(
\theta_{j,k}-\psi(\theta )_{j,k}\bigr) \bigl(
\theta_{j,k}^0-\psi(\theta)_{j,k}\bigr)
\leq2^{-7}\bigl\llVert \theta-\psi (\theta)\bigr\rrVert
_{L^2}^2.
\]
Since also $\sum_{j,k}\theta_{j,k}(\theta_{j,k}^0-\psi(\theta
)_{j,k})\leq\overline LQ\phi_n$, $\frac{1}2 -\llVert  \psi(\theta)\rrVert
_{L^2}\leq-1/4$ and
\[
-\sum_{j,k} \biggl(\frac{\theta_{j,k}}{\llVert  \theta\rrVert  _{L^2}}-
\frac{\psi
(\theta)_{j,k}}{\llVert  \psi(\theta)\rrVert  _{L^2}} \biggr)^2\leq-\frac{1}4\bigl\llVert
\theta-\psi(\theta)\bigr\rrVert _{L^2}^2+2 \biggl(
\frac{1}{\llVert  \theta\rrVert
_{L^2}}-\frac{1}{\llVert  \psi(\theta)\rrVert  _{L^2}} \biggr)^2,
\]
we obtain that $n^{-1}  (\mathcal L_n(\theta)-\mathcal L_n(\psi
(\theta)) )$ is less than
\begin{eqnarray*}
&&-\frac{3}{64}\bigl\llVert \theta-\psi(\theta)
\bigr\rrVert _{L^2}^2+\frac{1}2 \biggl(
\frac{1}{\llVert  \theta\rrVert  _{L^2}}-\frac{1}{\llVert  \psi(\theta)\rrVert  _{L^2}} \biggr)^2+2\overline LQ
\phi_n \biggl\llvert \frac{1}{\llVert  \theta\rrVert  _{L^2}}-\frac{1}{\llVert
\psi(\theta)\rrVert  _{L^2}}\biggr
\rrvert.
\end{eqnarray*}
Now, we need to distinguish two cases. First, assume that there is a
$(j,k)$ with $j\leq J_0$ and $\psi(\theta)_{j,k}\neq\theta_{j,k}$.
In this case, $\llVert  \theta-\psi(\theta)\rrVert  _{L^2}^2\geq Q^2\phi_n^2$. By
(\ref{eqthetaL^2ineqs}) and the choice of $Q$,
\begin{eqnarray*}
\frac{1}n \bigl(\mathcal L_n(\theta)-\mathcal
L_n\bigl(\psi(\theta)\bigr) \bigr) &\leq&-\frac{1}{64}\bigl
\llVert \theta-\psi(\theta)\bigr\rrVert _{L^2}^2+
\bigl(2^7\overline L^2+2^5\overline
L^2Q-2^{-5}Q^2\bigr)\phi_n^2
\\
&\leq& -\frac{1}{64}\bigl\llVert \theta-\psi(\theta)\bigr\rrVert
_{L^2}^2.
\end{eqnarray*}
Now, suppose the opposite, that is, whenever $\psi(\theta)_{j,k}\neq
\theta_{j,k}$ then $j>J_0$. By construction of $J_0$ [see (\ref
{eqdefofJ0})],
\begin{eqnarray*}
\biggl\llvert \frac{1}{\llVert  \theta\rrVert  _{L^2}}-\frac{1}{\llVert  \psi(\theta)\rrVert  _{L^2}}\biggr\rrvert &\leq& 4 \bigl
\llvert \bigl\langle\psi(\theta)+\theta, \psi(\theta)-\theta \bigr\rangle\bigr
\rrvert
\\
&\leq& 4 \max_{j>J_0, k\in I_j} \bigl\llvert \theta_{j,k}+\psi(
\theta )_{j,k}\bigr\rrvert \phi_n^{-1}\bigl\llVert
\theta-\psi(\theta)\bigr\rrVert _{L^2}^2
\\
&\leq& 2^{-7}(\overline L Q)^{-1} \phi_n^{-1}
\bigl\llVert \theta-\psi(\theta)\bigr\rrVert _{L^2}^2.
\end{eqnarray*}
Similar, we obtain
\[
\biggl(\frac{1}{\llVert  \theta\rrVert  _{L^2}}-\frac{1}{\llVert  \psi(\theta)\rrVert  _{L^2}} \biggr)^2
\leq2^{-5} \bigl\llVert \theta-\psi(\theta)\bigr\rrVert
_{L^2}^2,
\]
by using the Cauchy--Schwarz inequality instead. Therefore, in this
case, we also get $\frac{1}n  (\mathcal L_n(\theta)-\mathcal
L_n(\psi(\theta)) )\leq-\frac{1}{64}\llVert  \theta-\psi(\theta)\rrVert
_{L^2}^2$. This completes the proof.
\end{longlist}
\end{appendix}

\section*{Acknowledgements}
The comments and suggestions of two anonymous referees and an Associate
Editor that helped to improve considerably a previous version of this
work are gratefully acknowledged.



\printaddresses
\end{document}